\documentclass{article}

 \usepackage{amsmath,amssymb,amsthm,amsfonts,graphicx,color,epic}
 \usepackage[TABBOTCAP]{subfigure}
 \usepackage{epsfig}
 \usepackage{curves,eepic}
 \usepackage{bm}
 \usepackage{bbm}
 \usepackage{soul}
 \usepackage{rotating}

 \newtheorem{thm}{Theorem}

 \newtheorem{prop}{Proposition}
 
 \theoremstyle{definition}

 \newcommand{\E}{\mathbb{E}}
 \renewcommand{\P}{\mathbb{P}}

\begin{document}
 %TITLE

 \centerline{\Large Performance analysis of polling systems}
 \centerline{\Large with retrials and glue periods}

 \vspace{0.3cm}

% AUTHORS

 \centerline{Murtuza Ali Abidini$^*$, Onno Boxma$^*$, Bara Kim$^{**}$, Jeongsim Kim$^{***}$, Jacques
 Resing$^*$}

 \centerline{\small \it $^{*}$EURANDOM and Department of Mathematics and Computer
 Science}
 \centerline{\small \it Eindhoven University of Technology}
 \centerline{\small \it P.O. Box 513, 5600 MB Eindhoven, The Netherlands}
 \centerline{\small \it e-mail: m.a.abidini@tue.nl, o.j.boxma@tue.nl, j.a.c.resing@tue.nl}

 \vspace{0.2cm}

 \centerline{\small \it $^{**}$Department of Mathematics, Korea University}
 \centerline{\small \it 145 Anam-ro, Seongbuk-gu, Seoul, 02841, Korea}
 \centerline{\small \it e-mail: bara@korea.ac.kr}

 \vspace{0.2cm}

 \centerline{\small \it $^{***}$Department of Mathematics
 Education, Chungbuk National University}
 \centerline{\small \it 1 Chungdae-ro, Seowon-gu,
 Cheongju, Chungbuk, 28644, Korea}
 \centerline{\small \it e-mail: jeongsimkim@chungbuk.ac.kr}

\begin{abstract}
We consider gated polling systems with two special features: (i)
retrials, and (ii) glue or reservation periods. When a type-$i$
customer arrives, or retries, during a glue period of station $i$,
it will be served in the next visit period of the server to that
station. Customers arriving at station $i$ in any other period
join the orbit of that station and will retry after an
exponentially distributed time. Such polling systems can be used
to study the performance of certain switches in optical
communication systems.

For the case of exponentially distributed glue periods, we present
an algorithm to obtain the moments of the number of customers in
each station. For generally distributed glue periods, we consider
the distribution of the total workload in the system, using it to
derive a pseudo conservation law which in its turn is used to
obtain accurate approximations of the individual mean waiting
times. We also consider the problem of choosing the lengths of the
glue periods, under a constraint on the total glue period per
cycle, so as to minimize a weighted sum of the mean waiting times.
\end{abstract}


\numberwithin{equation}{section}

\section{Introduction}

This paper is devoted to the performance analysis of a class of
single server queueing systems with multiple customer types. Our
motivation is twofold: (i) to obtain insight into the performance
of certain switches in optical communication systems, and (ii) to
obtain insight into the effect of having particular reservation
periods, windows of opportunity during which a customer can make a
reservation for service. Our class of queueing systems combines
several features, viz., polling, retrials, and the new feature of
so-called glue periods or reservation periods. These will first be
discussed separately, while their relation to optical switching
will also be outlined.

Polling systems are queueing models in which a single server,
alternatingly, visits a finite number of, say, $N$ queues (or
stations) in some prescribed order. Polling systems have been
extensively studied in the literature. For example, various
different service disciplines (rules which describe the server's
behaviour while visiting a queue) have been considered, both for
models with and without switchover times between queues. We refer
to Takagi \cite{Takagi1,Takagi2} and Vishnevskii and Semenova
\cite{Vishnevskii} for literature reviews and to Boon, van der Mei
and Winands \cite{Boon}, Levy and Sidi \cite{Levy} and Takagi
\cite{Takagi3} for overviews of the applicability of polling
systems.

Switches in communication systems form an important application
area of polling systems. Here, packets must be routed from source
to destination, passing through a series of links and nodes. In
copper-based transmission links, packets from various sources are
time-multiplexed, and this may be modelled by a polling system. In
recent years optical networking has become very important, because
optical fibers offer major advantages with respect to copper
cables: huge bandwidth, ultra-low losses and an extra dimension,
viz., a choice of wavelengths.

When one wants to model  the performance of an optical switch by a
polling system \cite{Maier,Rogiest}, one is faced with the
following difficulty. Buffering of optical packets is not easy, as
photons can not wait. Whenever there is a need to buffer photons,
they are sent into a local fiber loop, thus providing a small
delay to the photons without losing or displacing them. If, at the
completion of the loop, a photon still needs to be buffered, it is
again sent into the fiber delay loop, etc. From a queueing
theoretic perspective, this raises the need to add the feature of
{\em retrial queue} to a polling system: instead of having a
queueing system with one server and $N$ ordinary queues, it has
one server and $N$ retrial queues. Retrial queues have received
much attention in the literature, see, e.g., the books by Falin
and Templeton \cite{Falin} and by Artalejo and Gomez-Corral
\cite{Artalejo}, but they have hardly been studied in the setting
of polling models. Langaris \cite{Langaris1,Langaris2,Langaris3}
has pioneered the study of polling models with retrial queues.
However, for our purpose - the performance analysis of optical
switches - his assumptions about the service discipline of the
server at the various queues are not suitable.

A third important feature in the present paper is that of
so-called glue or reservation periods. Just before the server
arrives at a station there is some glue period. Customers (both
new arrivals and retrying customers) arriving at the station
during this glue period ``stick" and will be served during the
visit of the server. Customers arriving  at the station in any
other period join the orbit of that station and will retry after
an exponentially distributed time. One motivation for studying
glue periods is the following. A sophisticated technology that one
might try to add to the use of fiber delay loops in optical
networking is {\em varying the speed of light by changing the
refractive index of the fiber loop}, cf. \cite{okawachi}. Using a
higher refractive index in a small part of the loop one can
achieve `slow light', which implies slowing down the packets. This
feature is in our model incorporated as glue periods, where we
slow down the packets arriving at the end of the fiber loop just
before the server arrives, so that they do not have to retry but
get served during the subsequent visit period. Not restricting
ourselves to optical networks, one can also interpret a glue
period as a reservation period, i.e., a period in which customers
can make a reservation at a station for service in the subsequent
visit period  of that station. In our model, the reservation
period immediately precedes the visit period, and could be seen as
the last part of a switchover period.

A first attempt to study a polling model which combines retrials
and glue periods is \cite{BoxmaResing}, which mainly focuses on
the case of a single server and a single station, but also
outlines how that analysis can be extended to the case of two
stations. In \cite{M_O_J_1} an $N$-station polling model with
retrials and with constant glue periods is considered, for the
case of gated service discipline at all stations. The gated
discipline is an important discipline in polling systems; it
implicates that the server, when visiting a station, serves
exactly those customers which were present upon his arrival. The
steady-state joint station size (i.e., the number of customers in
each station) distribution was derived in \cite{M_O_J_1}, both at
an arbitrary epoch and at beginnings of switchover, glue and visit
periods. In the current paper we present an algorithm to obtain
the moments of the station size for the case of exponentially
distributed glue periods. Using Little's law, that also gives mean
sojourn times. Thereafter for each individual station we allow
generally distributed glue periods and we focus our attention on
other performance measures next to station sizes. In particular,
we consider the steady-state distribution of the total workload in
the system, which leads us to a pseudo conservation law, i.e., an
exact expression for a weighted sum of the mean waiting times at
all stations. We use that pseudo conservation law to derive an
accurate approximation for the individual mean waiting times. We
further consider the problem of choosing the lengths of the glue
periods, given the total glue period in a cycle, so as to minimize
a weighted sum of the mean waiting times.

The rest of the paper is organized as follows. In Section
\ref{section_model} we present a model description. Section
\ref{section_old} contains a detailed analysis for generating
functions and moments of station sizes at different time epochs
when the glue periods are exponentially distributed. We also
present a numerical example, which in particular provides insight
into the behavior of the polling system in the case of long glue
periods. In Section \ref{section_new} we derive a pseudo
conservation law for a system with generally distributed glue
periods. Subsequently we use this pseudo conservation law for
deriving an approximation for the mean waiting times at all
stations. In its turn, this approximation is used to minimize
weighted sums of the mean waiting times by optimally choosing the
lengths of the glue periods, given the total glue period per
cycle. Finally, Section~\ref{section_concl} lists some topics for
further research.

 \section{The model} \label{section_model}

 We consider a single server cyclic polling system with retrials and
 so-called glue periods. This model was first introduced in \cite{BoxmaResing}
 for a single station vacation model
 and a two-station model with switchover times. Further in \cite {M_O_J_1} this was extended to an $N$-station
 model with switchover times. In both papers, the model was studied for deterministic glue periods.
 We index the stations by $i$, $i=1, \ldots, N$, in the
 order of server movement. For ease of presentation, all references to station indices greater
 than $N$ or less than $1$ are implicitly assumed to be modulo $N$.
 Customers arrive at station $i$ according to a Poisson
 process with rate $\lambda_i$; they are called type-$i$
 customers, $i=1,\ldots,N$.  The overall arrival rate is denoted by
 $\lambda=\lambda_1+\cdots+\lambda_N$. The service times at station $i$ are
 independent and identically distributed (i.i.d.) random variables
 with a generic random variable $B_i$, $i=1,\ldots,N$. Let
 $\tilde{B}_i(s)=\mathbb{E}[e^{-sB_i}]$ be the Laplace-Stieltjes transform
 (LST) of the service time distribution at station $i$.
 The switchover times from station $i$ to station $i+1$
 are i.i.d. random variables with a generic random variable $S_i$.
 Let $\tilde{S}_i(s)=\mathbb{E}[e^{-sS_i}]$ be the LST
 of the switchover time from station $i$ to station $i+1$, $i=1, \ldots, N$.
 The interarrival times, the service times, and the
 switchover times are assumed to be mutually independent.
 After a switch of the server to station $i$, there is a glue period
 for collecting retrying customers (which will be followed by the visit
 period of the server to station $i$).
 We assume that the successive glue periods at station $i$ are i.i.d. random variables
 with a generic random variable $G_i$. Let  $\tilde{G}_i(s)=\mathbb{E}[e^{-sG_i}]$ be the LST
 of the glue period distribution at station $i$.

 Each station consists of an orbit and a queue.
 When customers (both new arrivals and retrying customers) arrive at
 station $i$ during a glue period, they stick and wait in the queue to get served during the visit
 of the server to that station. When customers arrive at station $i$
 in any other period, they join the orbit of station $i$
 and will retry after a random amount of time. The inter-retrial time of each customer in
 the orbit of station $i$ is exponentially distributed with mean $\nu_i^{-1}$ and
 is independent of all other processes.

 A single server cyclically moves from one station to another serving the glued customers at each of the stations.
 The service discipline at all stations is gated. During the service
 period of station $i$, the server serves all glued customers in
 the queue of station $i$, i.e., all type-$i$ customers waiting at the end of
 the glue period (but none of those in orbit, and neither any new arrivals).

Let $(X_{1}^{(i)},X_{2}^{(i)},\ldots ,X_{N}^{(i)})$ denote the
vector of numbers of customers of type $1$ to type $N$ in the
system (hence in the orbit) at the start of a glue period of
station $i$, $i=1,\ldots,N$, in steady state. Further, let
$(Y_{1}^{(i)},Y_{2}^{(i)},\ldots ,Y_{N}^{(i)})$ denote the vector
of numbers of customers of type $1$ to type $N$ in the system at
the start of a visit period at station $i$, $i=1,\ldots,N$, in
steady state. We distinguish between those who are queueing
(glued) and those who are in the orbit of station $i$: We write
$Y_{i}^{(i)} = Y_{i}^{(iq)} + Y_{i}^{(io)}$, $i=1,\ldots ,N$,
where $q$ denotes in the queue and $o$ denotes in the orbit.
Finally, let $(Z_{1}^{(i)},Z_{2}^{(i)},\ldots ,Z_{N}^{(i)})$
denote the vector of numbers of customers of type $1$ to type $N$
in the system (hence in the orbit) at the start of a switchover
from station $i$ to station $i+1$, $i=1,\ldots ,N$, in steady
state.

 The utilization of the server at station $i$, $\rho_i$,
 is defined by $\rho_i = \lambda_i \mathbb{E}[B_i]$ and the total
 utilization of the server $\rho$ is given by $\rho=\sum_{i=1}^N \rho_i$.
 It can be shown that a necessary and sufficient condition for stability of
 this polling system is $\rho < 1$. We hence assume that $\rho<1$.

 The cycle length of station $i$, $i = 1,\ldots,N$  is defined as the time between two
 successive arrivals of the server at this station.
 The mean cycle length, $\mathbb{E}[C]$, is independent of the station
 involved (and the service discipline) and is given by
 \begin{eqnarray}
 \mathbb{E}[C]=\frac{\sum^N_{i=1} (\E[G_i]+\E[S_i])}{1-\rho},
 \label{eqn-C}
 \end{eqnarray}
 which can be derived as follows:
 Since the probability of the server being idle (in steady state) is
 $1-\rho$,
 and this equals $\frac{\sum^N_{i=1} (\E[G_i]+\E[S_i])}{\mathbb{E}[C]}$
 by the theory of regenerative processes, we have
 \eqref{eqn-C}.

\section{The polling system with retrials and exponential glue periods}
\label{section_old}

In \cite{M_O_J_1} the authors calculated the generating functions
and the mean values of the number of customers at different time
epochs when the glue periods are deterministic. In this section we
assume that the glue periods are \textit{exponentially}
distributed with mean $\E [G_i] =1/\gamma_i$, $i= 1, \ldots, N$.
We will derive a set of partial differential equations for the
joint generating function of the station size (i.e., the number of
customers in each station) and then obtain a system of linear
equations for the first and the second moments of the station
size. We also provide an iterative algorithm for solving the
system of linear equations.

Observe that the generating function for the vector of numbers of
arrivals at station $1$ to station $N$ during the service time of
a type-$i$ customer, $B_i$, is $\beta_i({\bf z}):=
\tilde{B}_i(\sum_{j=1}^{N}\lambda_j(1-z_j))$ for ${\bf z} =
(z_1,z_2, \ldots, z_N)$. Similarly, the generating function for
the vector of numbers of arrivals at
 station $1$ to station $N$ during a switchover time from station $i$ to station $i+1$,
$S_i$, is $\sigma_i({\bf z}):=\tilde{S}_i(\sum_{j=1}^{N}\lambda_j(1-z_j))$.

\subsection{Station size analysis at embedded time points}
In this subsection we study the steady-state joint distribution
and the mean of the numbers of customers in the system at the
start of a glue period, visit period and switchover period. Let us
define the following joint generating functions of the number of
customers in each station at the start of a glue period, visit
period and switchover period:
\begin{align*}
    \tilde{R}_{g}^{(i)}({\bf z}) &=\E [z_1^{X_{1}^{(i)}}
    z_2^{X_{2}^{(i)}}\cdots z_N^{X_{N}^{(i)}}], \\
    \tilde{R}_{v}^{(i)}({\bf z},w) &=\E
    [z_1^{Y_{1}^{(i)}}z_2^{Y_{2}^{(i)}}\cdots
    z_{i}^{Y_{i}^{(io)}}\cdots
    z_N^{Y_{N}^{(i)}}w^{Y_{i}^{(iq)}}], \\
    \tilde{R}_{s}^{(i)}({\bf z}) &=\E [z_{1}^{Z_{1}^{(i)}}
    z_2^{Z_{2}^{(i)}}\cdots z_N^{Z_{N}^{(i)}}],
\end{align*}
for ${\bf z} = (z_1,z_2, \ldots, z_N)$ with $|z_i| \leq 1$,
$i=1,\ldots,N$, and $|w| \leq 1$.

Let $M^o_i(t)$ represent the number of customers in the orbit of
station $i$, $i= 1, \cdots, N$ and $\Upsilon(t)$ the number of
glued customers, at time $t$. Further, let $\tau_j$ be the time at
which an arbitrary glue period starts at station $j$, $j= 1,
\ldots, N$. Note that $M_i^o(\tau_j)=X_i^{(j)}$. We define
\begin{align*}
    \phi_i({\bf z};w;t)&= \E[z_1^{M^o_1(\tau_i+t)}\cdots z_N^{M^o_N(\tau_i+t)} w^{\Upsilon(\tau_i+t)}\mathbbm{1}_{\{G_i>t\}}], \nonumber \\
    \phi_i({\bf z},w)&= \int^\infty_0 \phi_i({\bm z};w;t)dt.
\end{align*}
Then all the generating functions for the numbers of
customers in steady state described above, can be expressed in terms of
$\phi_{i}({\bf z},w)$, as shown below in Proposition \ref{sec5-prop1}.

\begin{prop} \label{sec5-prop1}
    The generating functions $\tilde{R}_{v}^{(i)}({\bf z},w), \tilde{R}_{s}^{(i)}({\bf
        z})$ and $\tilde{R}_{g}^{(i)}({\bf z})$ satisfy the following:
    \begin{align}
        \tilde{R}_{v}^{(i)}({\bf z},w) &=\gamma_i \phi_i({\bf
            z},w), \label{thm_1_b}\\
        \tilde{R}_{s}^{(i)}({\bf z}) &=\gamma_i \phi_i({\bf z},\beta_i({\bf
            z})), \label{thm_1_c}\\
        \tilde{R}_{g}^{(i)}({\bf z}) &=
        \gamma_{i-1} \sigma_{i-1}({\bf z})\phi_{i-1}({\bf z},\beta_{i-1}({\bf z})).
        \label{thm_1_new}
    \end{align}
\end{prop}

\noindent{\bf Proof.~}
Equation (\ref{thm_1_b}) is obtained as follows: By the law of total expectation,
\begin{align*}
    \tilde{R}_{v}^{(i)}({\bf z},w)
    &= \int^\infty_0 \E[z_1^{M^o_1(\tau_i+t)}\cdots z_N^{M^o_N(\tau_i+t)} w^{\Upsilon(\tau_i+t)}\mid G_i>t] \gamma_ie^{-\gamma_i t} dt\nonumber \\
    &= \int^\infty_0 \E[z_1^{M^o_1(\tau_i+t)}\cdots z_N^{M^o_N(\tau_i+t)} w^{\Upsilon(\tau_i+t)}\mathbbm{1}_{\{G_i>t\}}] \gamma_i dt \nonumber \\
    &= \gamma_i \phi_i({\bf z},w).
\end{align*}
To obtain (\ref{thm_1_c}), observe that the customers at the end
of a visit period are the customers in the orbit at the beginning
of that visit plus the customers who arrive during the service
times of the glued customers at the beginning of that visit. Hence
\begin{align*}
    \tilde{R}_{s}^{(i)}({\bf z})
    &= \E[z_1^{Y_{1}^{(i)}}z_2^{Y_{2}^{(i)}}\cdots z_i^{Y_{i}^{(io)}}\cdots z_N^{Y_{N}^{(i)}}[\beta_i({\bf z})]^{Y_{i}^{(iq)}}]\\
    &= \tilde{R}_v^{(i)} ({\bf z},\beta_i({\bf z}))\\
    &= \gamma_i \phi_i({\bf z},\beta_i({\bf z})).
\end{align*}
Also, to obtain (\ref{thm_1_new}), observe that the customers at
the end of a switchover from station $i-1$ to station $i$ are the customers
in the orbit at the beginning of that switchover plus the
customers who arrived during that switchover period. Hence
\begin{eqnarray*}
    \tilde{R}_{g}^{(i)}({\bf z})
    =\tilde{R}_s^{(i-1)}({\bf z}) \sigma_{i-1}({\bf z}),
\end{eqnarray*}
from which and (\ref{thm_1_c}) we get (\ref{thm_1_new}). \qed

We have the following result for the generating functions
$\phi_i({\bf z},w)$, $i=1,\ldots,N$.

\begin{thm} \label{thm1}
    The generating functions $\phi_i({\bf z},w)$, $i=1,\ldots,N$,
    satisfy the following equation:
    \begin{align}
        & \nu_i(w-z_i) \frac{\partial}{\partial z_i} \phi_i({\bf z},w)
        -\Big(\sum_{j=1,j\ne i}^{N}(\lambda_j(1-z_j)) +\lambda_i(1-w)+\gamma_i\Big)\phi_i({\bf
            z},w) \nonumber\\
        & \quad + \gamma_{i-1} \phi_{i-1}({\bf z},{\beta}_{i-1}({\bf z} ))
        \sigma_{i-1}({\bf z})=0.
        \label{thm_1}
    \end{align}
\end{thm}
{\noindent \bf Proof. }
Note that
\begin{align*}
    & \phi_i({\bf z};w;t+\Delta t) \\
    &= \E \big[z_1^{M^o_1(\tau_i+t+\Delta t)} \cdots z_N^{M^o_N(\tau_i+t+\Delta t)}
    w^{\Upsilon(\tau_i+t+\Delta t)}\mathbbm{1}_{\{G_i>t+\Delta t\}} \big] \nonumber \\
    &= \sum^\infty_{n_1=0} \cdots \sum^\infty_{n_N=0}\sum^\infty_{k=0}
    \P(M^o_1(\tau_i+t)=n_1,\ldots,M^o_N(\tau_i+t)=n_N,\Upsilon(\tau_i+t)=k,G_i>t) \\
    & \quad \times  \E \big[z_1^{M^o_1(\tau_i+t+\Delta t)}
    \cdots z_N^{M^o_N(\tau_i+t+\Delta t)} w^{\Upsilon(\tau_i+t+\Delta t)}
    \mathbbm{1}_{\{G_i>t+\Delta t\}} \big| M^o_1(\tau_i+t)=n_1,\ldots, \\
    &\qquad \qquad M^o_N(\tau_i+t)=n_N, \Upsilon(\tau_i+t)=k,G_i>t \big] \\
    &= \sum^\infty_{n_1=0} \cdots \sum^\infty_{n_N=0}\sum^\infty_{k=0}
    \P(M^o_1(\tau_i+t)=n_1,\ldots,M^o_N(\tau_i+t)=n_N,\Upsilon(\tau_i+t)=k,G_i>t) \\
    & \quad \times  z_1^{n_1}\cdots z_{i-1}^{n_{i-1}}z_{i+1}^{n_{i+1}}\cdots z_N^{n_N} w^k
    ((1-e^{-\nu_i\Delta t})w+ e^{-\nu_i \Delta t} z_{i})^{n_i}
    e^{-(\sum_{j=1,j\ne i}^{N}(\lambda_j(1-z_j)) +\lambda_i(1-w))\Delta t} e^{-\gamma_i \Delta t} \\
    &=  e^{-\big(\sum_{j=1,j\ne i}^{N}(\lambda_j(1-z_j)) +\lambda_i(1-w)+\gamma_i\big)\Delta t}
    \phi_i(z_1, \ldots, z_{i-1}, z_i+(1-e^{-\nu_i\Delta t})(w-z_i),z_{i+1},\ldots,z_N;w;t).
\end{align*}
Thus, we have
$$\frac{\partial}{\partial t} \phi_i({\bf z};w;t)
= \nu_i(w-z_i) \frac{\partial}{\partial z_i} \phi_i({\bf z};w;t)
-\Big(\sum_{j=1,j\ne i}^{N}(\lambda_j(1-z_j))
+\lambda_i(1-w)+\gamma_i\Big)\phi_i({\bf z};w;t).$$ Since
$\phi_i({\bf z};w;0)= \E [z_1^{X_{1}^{(i)}}
z_2^{X_{2}^{(i)}}\cdots z_N^{X_{N}^{(i)}}] =
\tilde{R}_{g}^{(i)}({\bf z})= \gamma_{i-1} \sigma_{i-1}({\bf z})
\phi_{i-1}({\bf z},\beta_{i-1}({\bf z}))$ and $\phi_i({\bf
z};w;\infty)=0$, integrating the above equation with respect to
$t$ from 0 to $\infty$ yields
\begin{align*}
    -\gamma_{i-1} \sigma_{i-1}({\bf z}) \phi_{i-1}({\bf z},\beta_{i-1}({\bf z}))=& \nu_i(w-z_i)
    \frac{\partial}{\partial z_i} \phi_i({\bf z},w) \\
    &-\Big(\sum_{j=1,j\ne i}^{N}\lambda_j(1-z_j) +\lambda_i(1-w)+\gamma_i\Big)\phi_i({\bf z},w).
\end{align*}
This completes the proof. \qed

We now calculate the mean value of the station sizes at embedded
time points using the differential equation \eqref{thm_1}. For an
$N$-tuple $\bm l=(l_1, \ldots, l_N)$ of nonnegative integers, we
define
\begin{align*}
    |\bm l| = l_1+\cdots+l_N,~~ {\bm l}! = l_1!l_2!\cdots l_N!,
\end{align*}
and ${\bf z}^{\bm l}=z_1^{l_1}z_2^{l_2}\cdots z_N^{l_N}$. With
this notation, we define the following scaled moment:
\begin{align*}
    \Phi_i^{(\bm l,m)}&= \frac{1}{{\bm l}!m!}\frac{\partial^{|\bm l|+m}}
    {\partial {\bf z}^{\bm l}\partial
        w^m}\phi_i({\bf z}, w) \Big|_{{\bf z}={\bm 1}-, {w}=1-},
\end{align*}
where $\partial {\bf z}^{\bm l}=\partial z_1^{l_1} \cdots \partial
z_N^{l_N}$, and ${\bm 1}$ is the $N$-dimensional row vector with
all its components equal to one. The first scaled moments of
$\phi_i({\bf z},w)$, $i=1,2,\ldots,N$, can be obtained from the
following theorem.

\begin{thm} \label{sec5-new-thm2}
    We have
    \begin{itemize}
        \item[(i)] $\Phi_i^{({\bm 0},0)}=\frac{1}{\gamma_i}$, $i=1,\ldots,N$.
        \item[(ii)] $\Phi_{i}^{({\bm 1}_j,0)}$ and $\Phi_{i}^{({\bm
                0},1)}$, $0\le i,j \le N$, are given by the following recursion:
        for $j=1,\ldots, N$,
        \begin{align}
            \Phi_{j}^{({\bm 0},1)} &= \frac{\lambda_j}{\gamma_j}\E[C], \label{sec3-eqn1-1}\\
            \Phi_j^{({\bm 1}_j,0)} &= \frac{\lambda_j}{\nu_j}\Big(\E[C]-\frac{1}{\gamma_j}\Big),
            \label{sec3-eqn1-2}\\
            \Phi_{i}^{({\bm 1}_j,0)} &= \frac{\gamma_{i+1}}{\gamma_{i}}\Phi_{i+1}^{({\bm
                    1}_j,0)}+ \frac{\lambda_j}{\gamma_{i}}
            \Big((\delta_{i,j-1}-\rho_{i})\E[C]-\frac{1}{\gamma_{i+1}}-\mathbb{E}[S_{i}]\Big),
            ~~i=j-1,j-2,\ldots,j-N+1, \label{sec3-eqn1-3}
        \end{align}
        where ${\bm 0}$ is the $N$-dimensional row vector with all its
        elements equal to zero, ${\bm 1}_j$ is the $N$-dimensional row vector whose
        $j$th element is one and all other elements are zero,
        and $\delta_{ij}$ is the Kronecker delta.
        Note that if $i$ is nonpositive in (\ref{sec3-eqn1-3}), then it
        is interpreted as $i+N$.
    \end{itemize}
\end{thm}

\noindent{\bf Proof.~} Taking the partial derivative of Equation
\eqref{thm_1} with respect to $z_j$ and putting ${\bf z}={\bm 1}-,
w= 1-$, we have
\begin{align}
    -\nu_i \delta_{ij} \Phi_i^{({\bm
            1}_i,0)}+\frac{(1-\delta_{ij})\lambda_j}{\gamma_i}-\gamma_i\Phi_i^{({\bm
            1}_j,0)}+\gamma_{i-1}\Phi_{i-1}^{({\bm
            1}_j,0)}+\gamma_{i-1}\lambda_j\mathbb{E}[B_{i-1}]\Phi_{i-1}^{({\bm
            0},1)}+\lambda_j\mathbb{E}[S_{i-1}]=0. \label{mom-eqn1}
\end{align}
Taking the partial derivative of Equation \eqref{thm_1} with
respect to $w$ and putting ${\bf z}={\bm 1}-, w=1-$ yields
\begin{align}
    \nu_i \Phi_i^{({\bm 1}_i,0)}+\frac{\lambda_i}{\gamma_i}-\gamma_i\Phi_{i}^{({\bm
            0},1)}=0. \label{mom-eqn2}
\end{align}
Summing (\ref{mom-eqn1}) over $i=1, \ldots,N$, we have
\begin{align}
    -\nu_j \Phi_j^{({\bm 1}_j,0)}+\lambda_j\sum_{i\neq
        j}\frac{1}{\gamma_i}+\lambda_j \sum_{i=1}^N \gamma_i\mathbb{E}[B_i]\Phi_{i}^{({\bm
            0},1)}+\lambda_j \sum_{i=1}^N \mathbb{E}[S_i]=0. \label{mom-eqn3}
\end{align}
Adding (\ref{mom-eqn2}) and (\ref{mom-eqn3}) and multiplying the resulting equation by
$\mathbb{E}[B_j]$ yields
\begin{align*}
    \rho_j \sum_{i=1}^N\Big(\frac{1}{\gamma_i}+\mathbb{E}[S_i]\Big)
    -\gamma_j\mathbb{E}[B_j]\Phi_{j}^{({\bm 0},1)}
    +\rho_j\sum_{i=1}^N \gamma_i\mathbb{E}[B_i]\Phi_{i}^{({\bm 0},1)}=0,
\end{align*}
and summing this over $j=1,\ldots,N$ gives
\begin{align}
    \sum_{i=1}^N \gamma_i\mathbb{E}[B_i]\Phi_{i}^{({\bm 0},1)}=\rho \E[C], \label{mom-eqn4}
\end{align}
where we have used (\ref{eqn-C}).
Plugging (\ref{mom-eqn4}) into (\ref{mom-eqn3}) leads to
\begin{align*}
    \Phi_j^{({\bm
    1}_j,0)}=\frac{\lambda_j}{\nu_j}\Big(\E[C]-\frac{1}{\gamma_j}\Big),
    \end{align*}
    which is (\ref{sec3-eqn1-2}).
Inserting this equation into (\ref{mom-eqn2}) yields
(\ref{sec3-eqn1-1}). When $i=j$ in Equation (\ref{mom-eqn1}), we
have
\begin{align}
    \Phi_{j-1}^{({\bm 1}_j,0)}=\frac{\gamma_j}{\gamma_{j-1}}\Phi_j^{({\bm
            1}_j,0)}+ \frac{\lambda_j}{\gamma_{j-1}}
    \Big((1-\rho_{j-1})\E[C]-\frac{1}{\gamma_j}-\mathbb{E}[S_{j-1}]\Big).  \label{mom-eqn7}
\end{align}
On the other hand, when $i\neq j$, i.e., $i=j-1,j-2,
\ldots,j-N+1$, in Equation (\ref{mom-eqn1}), we have
\begin{align}
    \Phi_{i-1}^{({\bm 1}_j,0)}=\frac{\gamma_i}{\gamma_{i-1}}\Phi_i^{({\bm
            1}_j,0)}+ \frac{\lambda_j}{\gamma_{i-1}}
    \Big(-\rho_{i-1}\E[C]-\frac{1}{\gamma_i}-\mathbb{E}[S_{i-1}]\Big). \label{mom-eqn8}
\end{align}
Finally, (\ref{sec3-eqn1-3}) follows from (\ref{mom-eqn7}) and
(\ref{mom-eqn8}). \qed

\vspace{0.4cm}
Next, we calculate $\Phi_i^{({\bm l},m)}$ for $|\bm l|+m \ge 2$.
Equation \eqref{thm_1} can be written as
\begin{align*}
    & \big(\nu_i(w-1) -\nu_i(z_i-1)\big)\frac{\partial}{\partial z_i} \phi_i({\bf z},w)
    +\Big(\lambda_i(w-1)+\sum_{j=1,j\neq i}^N\lambda_j(z_j-1)-\gamma_i\Big)\phi_i({\bf
        z},w) \\
    &+ \gamma_{i-1} \phi_{i-1}({\bf z},\beta_{i-1}({\bf z}))
    \sigma_{i-1}({\bf z})=0.
\end{align*}
From this we get
\begin{align}
    (\gamma_i+l_i\nu_i)\Phi_i^{({\bm l},m)} =& \mathbbm{1}_{\{m \ge
        1\}}(l_i+1)\nu_i \Phi_i^{({\bm l}+{\bm 1}_i,m-1)}
    +\mathbbm{1}_{\{m \ge 1\}}\lambda_i \Phi_i^{({\bm l},m-1)} \nonumber\\
    &+\sum_{j\neq i}\mathbbm{1}_{\{l_j \ge 1\}}\lambda_j \Phi_i^{({\bm l}-{\bm
            1}_j,m)}+\mathbbm{1}_{\{m=0\}}\gamma_{i-1}\sum_{{\bm l}'\le {\bm
            l}}\sum_{k=0}^{|{\bm l}-{\bm l}'|} \Phi_{i-1}^{({\bm l}',k)}
    \Gamma_{i-1,k}^{({\bm l}-{\bm l}')}, \label{sec3-eqn-new1}
\end{align}
where $\Gamma_{i,m}^{(\bm l)}=\frac{1}{{\bm l}!}\frac{\partial^{|{\bm l}|}}{\partial{\bf z}^{\bm l}}
\big((\beta_i({\bf z})-1)^m\sigma_i({\bf z})\big)\big|_{{\bf z}={\bm 1}-}$
and the inequality ${\bm l}'\le {\bm l}$ is interpreted componentwise.
Therefore, from (\ref{sec3-eqn-new1}) we have the following proposition.
\begin{prop}
    For $|\bm l|+m \ge 2$,
    \begin{align}
        \Phi_i^{({\bm l},m)} =&
        \frac{\mathbbm{1}_{\{m \ge 1\}}\lambda_i}{\gamma_i+l_i\nu_i} \Phi_i^{({\bm l},m-1)}
        +\sum_{j\neq i}\frac{\mathbbm{1}_{\{l_j \ge 1\}}\lambda_j}{\gamma_i+l_i\nu_i}
        \Phi_i^{({\bm l}-{\bm 1}_j,m)}
        +\frac{\mathbbm{1}_{\{m=0\}}\gamma_{i-1}}{\gamma_i+l_i\nu_i}\sum_{{\bm l}'\le {\bm
                l}}\sum_{k=0}^{|{\bm l}-{\bm l}'|-1} \Phi_{i-1}^{({\bm l}',k)}
        \Gamma_{i-1,k}^{({\bm l}-{\bm l}')} \nonumber\\
        &+ \frac{\mathbbm{1}_{\{m \ge 1\}}(l_i+1)\nu_i}{\gamma_i+l_i\nu_i} \Phi_i^{({\bm l}+{\bm 1}_i,m-1)}
        + \frac{\mathbbm{1}_{\{m=0\}}\gamma_{i-1}}{\gamma_i+l_i\nu_i}\sum_{{\bm l}'\le {\bm
                l}}\Phi_{i-1}^{({\bm l}',|{\bm l}-{\bm l}'|)}
        \Gamma_{i-1,|{\bm l}-{\bm l}'|}^{({\bm l}-{\bm l}')}. \label{sec3-eqn-new2}
    \end{align}
\end{prop}

We note that (\ref{sec3-eqn-new2}) is a system of linear equations
for $\Phi_i^{({\bm l},m)}$. This system of linear equations can be
solved by the Gaussian elimination method. However, we will use an
iterative method to solve the system of linear equations
(\ref{sec3-eqn-new2}). In the following theorem the iterative
algorithm is presented and the convergence of iteration is
guaranteed.

\begin{thm} \label{sec5-1-iteration}
    For $i, {\bm l}, m,n$ with $i=1,\ldots, N$, $|\bm l|+m=k$,
    $n=0,1,\ldots$, define $\Phi_i^{({\bm l},m)}(n)$ as follows:
    \begin{align*}
        \Phi_i^{({\bm l},m)}(0) =& 0, \\
        \Phi_i^{({\bm l},m)}(n) =&
        \frac{\mathbbm{1}_{\{m \ge 1\}}\lambda_i}{\gamma_i+l_i\nu_i} \Phi_i^{({\bm
                l},m-1)}
        +\sum_{j\neq i}\frac{\mathbbm{1}_{\{l_j \ge 1\}}\lambda_j}{\gamma_i+l_i\nu_i}
        \Phi_i^{({\bm l}-{\bm 1}_j,m)} \\
        &+\frac{\mathbbm{1}_{\{m=0\}}\gamma_{i-1}}{\gamma_i+l_i\nu_i}\sum_{{\bm l}'\le {\bm
                l}}\sum_{k=0}^{|{\bm l}-{\bm l}'|-1} \Phi_{i-1}^{({\bm l}',k)}
        \Gamma_{i-1,k}^{({\bm l}-{\bm l}')}
        + \frac{\mathbbm{1}_{\{m \ge 1\}}(l_i+1)\nu_i}{\gamma_i+l_i\nu_i} \Phi_i^{({\bm l}+{\bm
                1}_i,m-1)}(n-1) \\
        &+ \frac{\mathbbm{1}_{\{m=0\}}\gamma_{i-1}}{\gamma_i+l_i\nu_i}\sum_{{\bm l}'\le {\bm
                l}}\Phi_{i-1}^{({\bm l}',|{\bm l}-{\bm l}'|)}(n-1)
        \Gamma_{i-1,|{\bm l}-{\bm l}'|}^{({\bm l}-{\bm l}')}, ~~~n \ge 1.
    \end{align*}
    Then we have that
    \begin{itemize}
        \item[(i)] $\Phi_i^{({\bm l},m)}(n)$ is nondecreasing in $n$.
        \item[(ii)] $\lim_{n \to \infty}\Phi_i^{({\bm l},m)}(n)=\Phi_i^{({\bm
                l},m)}$.
    \end{itemize}
\end{thm}

\noindent{\bf Proof.~} By induction on $n$, we have that $\Phi_i^{({\bm
        l},m)}(n)$ is increasing in $n$ and $\Phi_i^{({\bm l},m)}(n) \le \Phi_i^{({\bm
        l},m)}$ for all $n$. Thus (i) is proved. Moreover, $\lim_{n \to \infty}\Phi_i^{({\bm
        l},m)}(n)$ exists and $\lim_{n \to \infty}\Phi_i^{({\bm l},m)}(n) \le \Phi_i^{({\bm
        l},m)}$.
Suppose that $\{\lim_{n\to \infty}\Phi_i^{({\bm l},m)}(n):
i=1,\ldots,N, |\bm l|+m=k\}$ and $\{\Phi_i^{({\bm l},m)}:
i=1,\ldots,N, |\bm l|+m=k\}$ are different solutions of the system
of equations (\ref{sec3-eqn-new2}). Then $\{a\lim_{n \to
\infty}\Phi_i^{({\bm l},m)}(n) +(1-a)\Phi_i^{({\bm l},m)}:
i=1,\ldots,N, |\bm l|+m=k \}$ is a solution for any $a \in
\mathbb{R}$. Since $\lim_{n \to \infty}\Phi_i^{({\bm l},m)}(n) \le
\Phi_i^{({\bm l},m)}$, there exists $a$ such that
\begin{align*}
    a\lim_{n \to \infty}\Phi_i^{({\bm l},m)}(n) +(1-a) \Phi_i^{({\bm
            l},m)} \ge 0
\end{align*}
for all $i=1,\ldots,N$ and $({\bm l}, m)$ with $|\bm l|+m=k$ and
\begin{align*}
    a\lim_{n \to \infty}\Phi_i^{({\bm l},m)}(n) +(1-a) \Phi_i^{({\bm
            l},m)}= 0
\end{align*}
for some $i=1,\ldots,N$ and $({\bm l}, m)$ with $|\bm l|+m=k$.
Hence there exists a nonnegative (vector) solution of
(\ref{sec3-eqn-new2}) with a zero component, which is a
contradiction, because $ \frac{\mathbbm{1}_{\{m \ge
1\}}\lambda_i}{\gamma_i+l_i\nu_i} \Phi_i^{({\bm
        l},m-1)}
+\sum_{j\neq i}\frac{\mathbbm{1}_{\{l_j \ge 1\}}\lambda_j}{\gamma_i+l_i\nu_i}
\Phi_i^{({\bm l}-{\bm 1}_j,m)}
+\frac{\mathbbm{1}_{\{m=0\}}\gamma_{i-1}}{\gamma_i+l_i\nu_i}\sum_{{\bm l}'\le {\bm
        l}}\sum_{k=0}^{|{\bm l}-{\bm l}'|-1} \Phi_{i-1}^{({\bm l}',k)}
\Gamma_{i-1,k}^{({\bm l}-{\bm l}')}$ is positive for all
$i=1,\cdots,N$ and $({\bm l}, m)$ with $|\bm l|+m=k$. Therefore,
$\lim_{n \to \infty}\Phi_i^{({\bm l},m)}(n)=\Phi_i^{({\bm l},m)}$
for all $i=1,\ldots,N$ and $({\bm l}, m)$ with $|\bm l|+m=k$. \qed

\subsection{Station size analysis at arbitrary time points}

In the previous subsection we have found the generating functions
of the number of customers at the beginning of glue periods, visit
periods, and switchover periods in terms of $\phi_i({\bf z},w)$.
We now represent the generating function of the number of
customers at arbitrary time points in terms of $\phi_i({\bf
z},w)$, as shown below in Theorem \ref{thm2}. This will allow us
to obtain the moments of the station size distribution at
arbitrary time points.

\begin{thm} \label{thm2}
    \begin{itemize}
        \item[(a)] The joint generating function, $R^{(i)}_{s}({\bf z})$, of the number of
        customers in the orbit at an arbitrary time point in a
        switchover period from station $i$ is given by
        \begin{equation}
            R^{(i)}_{s}({{\bf z}}) =  \frac{\gamma_i}{\mathbb{E}[S_i]}\phi_i({\bf z},\beta_i({\bf
                z}))\frac{1-\sigma_i({\bf z})}{\sum_{j=1}^N\lambda_j(1-z_{j})}.
            \label{sec5-2-neweqn3}
        \end{equation}
        \item[(b)] The joint generating function, $R^{(i)}_{g}({\bf z},w)$, of the number of
        customers in the queue and in the orbit at an arbitrary time point in a
        glue period of station $i$ is given by
        \begin{equation}
            % \hspace{-0.3in}
            R^{(i)}_{g}({\bf z},w) = \gamma_i \phi_i({\bf z},w). \label{sec5-2-neweqn4}
        \end{equation}
        \item[(c)] The joint generating function, $R^{(i)}_{v}({\bf z},w)$, of the number of
        customers in the queue and in the orbit at an arbitrary time point in a
        visit period of station $i$ is given by
        \begin{eqnarray}
            R^{(i)}_{v}({\bf{z}},w) &=&\frac{\gamma_i}{\rho_i\mathbb{E}[C]}
            \frac{\phi_i({\bf z},w) -\phi_i({\bf z},\beta_i({\bf z}))}{
                w - \beta_i({\bf z})} \frac{1-\beta_i ({\bf z})}
            {\sum_{j=1}^N\lambda_j(1-z_{j})}. \label{sec5-2-neweqn5}
        \end{eqnarray}
    \end{itemize}
\end{thm}

{\noindent \bf Proof. } (a) Notice that the number of customers in
the orbit at an arbitrary time point in a switchover period from
station $i$ is the sum of two independent terms: the number of
customers at the beginning of the switchover period and the number
of customers who arrived during the elapsed switchover period. The
generating function of the former is $\tilde{R}_{s}^{(i)}({\bf
z})$ and the generating function of the latter is given by
$\frac{1-\sigma_i({\bf z})}{\E
[S_i]\big(\sum_{j=1}^N\lambda_j(1-z_{j})\big)}$. Thus
\begin{align*}
    R^{(i)}_{s}({{\bf z}})=\tilde{R}_{s}^{(i)}({\bf z})\frac{1-\sigma_i({\bf z})} {\E
        [S_i]\big(\sum_{j=1}^N\lambda_j(1-z_{j})\big)},
\end{align*}
from which and (\ref{thm_1_c}) we get (\ref{sec5-2-neweqn3}).

(b) By the theory of Markov regenerative processes,
\begin{align*}
    R^{(i)}_{g}({\bf z},w)=\gamma_i\int_0^\infty \phi_i({\bf
        z};w;t)dt.
\end{align*}
This yields (\ref{sec5-2-neweqn4}).

(c) Notice that the number of customers in the system at an
arbitrary time point in a visit period consists of two parts: the
number of customers in the system at the beginning of the service
of the customer currently in service and the number of customers
who arrived during the elapsed time of the current service. The
generating function of the former is given by (see Remark 3 of
\cite{M_O_J_1} for a detailed proof)
\begin{align}
    \frac{\tilde{R}_{v}^{(i)}({\bf z},w) -\tilde{R}_{v}^{(i)}({\bf z},\beta_i({\bf z}))}
    {\mathbb{E}[Y_i^{(iq)}](w - \beta_i({\bf z}))}=
    \frac{\gamma_i}{\mathbb{E}[Y_i^{(iq)}]}\frac{\phi_i({\bf z},w)
        -\phi_i({\bf z},\beta_i({\bf z}))}{w - \beta_i({\bf z})},
    \label{sec5-2-new-eqn1}
\end{align}
and the generating function of the latter is given by
\begin{align}
    \frac{1-\beta_i ({\bf z})}
    {\E[B_i]\left(\sum_{j=1}^N\lambda_j(1-z_{j})\right)}. \label{sec5-2-new-eqn2}
\end{align}
From (\ref{sec5-2-new-eqn1}) and (\ref{sec5-2-new-eqn2}) we have
\begin{align*}
    R^{(i)}_{v}({\bf{z}},w)=\frac{\gamma_i}{\mathbb{E}[Y_i^{(iq)}]\E[B_i]}
    \frac{\phi_i({\bf z},w) -\phi_i({\bf z},\beta_i({\bf z}))}{w - \beta_i({\bf z})}
    \frac{1-\beta_i ({\bf z})} {\sum_{j=1}^N\lambda_j(1-z_{j})}.
\end{align*}
Since
$\rho_i=\frac{\mathbb{E}[Y_i^{(iq)}]\E[B_i]}{\mathbb{E}[C]}$,
(\ref{sec5-2-neweqn5}) follows from the above equation. \qed

We introduce the following scaled moments:
\begin{align*}
    \Psi_{g,i}^{(\bm l,m)} &= \frac{1}{{\bm l}!m!}\frac{\partial^{|\bm l|+m}}
    {\partial {\bf z}^{\bm l}\partial
        w^m} R^{(i)}_{g}({\bf z},w)\Big|_{{\bf z}={\bm 1}-, {w}=1-}, \\
    \Psi_{v,i}^{(\bm l,m)} &= \frac{1}{{\bm l}!m!}\frac{\partial^{|\bm l|+m}}
    {\partial {\bf z}^{\bm l}\partial
        w^m} R^{(i)}_{v}({\bf z},w)\Big|_{{\bf z}={\bm 1}-, {w}=1-}, \\
    \Psi_{s,i}^{(\bm l)} &= \frac{1}{{\bm l}!}\frac{\partial^{|\bm l|}}
    {\partial {\bf z}^{\bm l}} R^{(i)}_{s}({{\bf z}})\Big|_{{\bf z}={\bm 1}-}.
\end{align*}
These moments satisfy the following theorem, which can be derived by
using Equations (\ref{sec5-2-neweqn3}), (\ref{sec5-2-neweqn4}) and
(\ref{sec5-2-neweqn5}).

\begin{thm} \label{sec5-new-thm1}
    We have
    \begin{itemize}
        \item[(i)] $\Psi_{g,i}^{(\bm l,m)}=\gamma_i \Phi_i^{(\bm
        l,m)}$.

        \item[(ii)] $\Psi_{v,i}^{(\bm l,m)}=\frac{\gamma_i}{\rho_i\E[C]} \sum_{{\bm l}'\le {\bm
                l}}\sum_{k=0}^{|{\bm l}-{\bm l}'|}\Phi_{i}^{({\bm l}',m+k+1)}
        \eta_{i,k}^{({\bm l}-{\bm l}')}$,
        where $\eta_{i,m}^{(\bm l)}=\frac{1}{{\bm l}!}\frac{\partial^{|{\bm l}|}}{\partial{\bf z}^{{\bm l}}}
        \Big(\frac{-(\beta_i ({\bf z})-1)^{m+1}}{\sum_{j=1}^N\lambda_j(1-z_{j})}\Big)\Big|_{{\bf z}={\bf
                1}-}$.

        \item[(iii)] $\Psi_{s,i}^{(\bm l)}=\frac{\gamma_i}{\E[S_i]} \sum_{{\bm l}'\le {\bm
                l}}\Theta_i^{({\bm l}')} \zeta_{i}^{({\bm l}-{\bm l}')}$,
        where $\zeta_{i}^{(\bm l)}=\frac{1}{{\bm l}!}\frac{\partial^{|{\bm l}|}}{\partial{\bf z}^{{\bm l}}}
        \Big(\frac{1-\sigma_i ({\bf z})}{\sum_{j=1}^N\lambda_j(1-z_{j})}\Big)\Big|_{{\bf z}={\bf
                1}-}$
        and $\Theta_{i}^{(\bm l)}=\frac{1}{{\bm l}!}\frac{\partial^{|{\bm l}|}}{\partial{\bf z}^{{\bm l}}}
        \phi_i({\bf z},\beta_i({\bf z}))\big|_{{\bf z}={\bf 1}-}$.
        Moreover, $\Theta_{i}^{(\bm l)}$ is given by
        \begin{align*}
            \Theta_i^{(\bm l)} = \sum_{{\bm l}'\le {\bm
                    l}}\sum_{k=0}^{|{\bm l}-{\bm l}'|}\Phi_{i}^{({\bm l}',k)}
            \Delta_{i,k}^{({\bm l}-{\bm l}')},
        \end{align*}
        where $\Delta_{i,m}^{(\bm l)}=\frac{1}{{\bm l}!}\frac{\partial^{|{\bm l}|}}{\partial{\bf z}^{{\bm l}}}
        (\beta_i({\bf z})-1)^m\big|_{{\bf z}={\bf 1}-}$.
    \end{itemize}
\end{thm}

From now on we obtain the first and second moments of the station
sizes of each type of customers in steady state. Let $M_i^{o}$ and
$\Upsilon$ be the steady state random variables corresponding to
$M_i^{o}(t)$ and $\Upsilon(t)$, respectively. That is, $M_i^{o}$
is the number of customers in the orbit of station $i$ in steady
state and $\Upsilon$ is the number of glued customers in steady
state. Let $M_i^{oq}$ be the number of customers in the orbit of
station $i$ plus the glued customers in the queue of station $i$
in steady state, and $M_i$ be the number of customers in station
$i$ (including the customer in service at station $i$) in steady
state. Moreover, we define the following indicator random
variables: for $i=1, \ldots,N$,
\begin{align*}
    I_{v,i} &=\left\{%
    \begin{array}{ll}
        1 & \hbox{ if the server is serving at station $i$ in steady state,} \\
        0 & \hbox{ otherwise,} \\
    \end{array}%
    \right. \\
    I_{g,i} &=\left\{%
    \begin{array}{ll}
        1 & \hbox{ if the server is in the glue period of station $i$ in steady state,} \\
        0 & \hbox{ otherwise,} \\
    \end{array}%
    \right. \\
    I_{s,i} &=\left\{%
    \begin{array}{ll}
        1 & \hbox{ if the server is switching from station $i$ to station $i+1$ in steady state,} \\
        0 & \hbox{ otherwise.} \\
    \end{array}%
    \right.
\end{align*}
Then we have that for $i=1,\ldots,N$,
\begin{align*}
    M_i^o &= \sum_{k=1}^N M_i^o (I_{v,k}+I_{g,k}+I_{s,k}), \\
    M_i^{oq} &=M_i^o+\Upsilon (I_{v,i}+I_{g,i}), \\
    M_i &= M_i^{oq}+I_{v,i}.
\end{align*}
Therefore, the mean station sizes, $\mathbb{E}[M_i^o],
\mathbb{E}[M_i^{oq}]$, and $\mathbb{E}[M_i]$, $i=1,\ldots,N$, are
given by
\begin{align}
    \mathbb{E}[M_i^o] &= \sum_{k=1}^N \Big(\rho_k\Psi_{v,k}^{({\bm 1}_i,0)}
    + \frac{\mathbb{E}[G_k]}{\mathbb{E}[C]}\Psi_{g,k}^{({\bm 1}_i,0)}
    +\frac{\mathbb{E}[S_k]}{\mathbb{E}[C]}\Psi_{s,k}^{({\bm 1}_i)}\Big), \label{sec5-2-mean-eqn1}\\
    \mathbb{E}[M_i^{oq}] &= \mathbb{E}[M_i^o]+\rho_i\Psi_{v,i}^{({\bm 0},1)}
    + \frac{\mathbb{E}[G_i]}{\mathbb{E}[C]}\Psi_{g,i}^{({\bm 0},1)}, \label{sec5-2-mean-eqn2}\\
    \mathbb{E}[M_i] &= \mathbb{E}[M_i^{oq}]+\rho_i.
    \label{sec5-2-mean-eqn3}
\end{align}

Now, in order to obtain the second moments of the station sizes,
$\mathbb{E}[M_i^{o}M_j^{o}], \mathbb{E}[M_i^{oq}M_j^{oq}]$, and
$\mathbb{E}[M_iM_j]$, $i,j=1,\ldots,N$, note that
\begin{align*}
    M_i^{o}M_j^{o} &=\sum_{k=1}^N
    M_i^{o}M_j^{o}(I_{v,k}+I_{g,k}+I_{s,k}),\\
    M_i^{oq}M_j^{oq} &=
    M_i^{o}M_j^{o}+M_i^o \Upsilon(I_{v,j}+I_{g,j})+M_j^o\Upsilon(I_{v,i}+I_{g,i})
    +\Upsilon^2(I_{v,i}+I_{g,i})\mathbbm{1}_{\{i=j\}}, \\
    M_iM_j &= M_i^{oq}M_j^{oq} +M_i^{oq}I_{v,j}+M_j^{oq}I_{v,i}
    +I_{v,i}\mathbbm{1}_{\{i=j\}}.
\end{align*}
Therefore, the second moments of the station sizes are given by
\begin{align}
    \mathbb{E}[M_i^{o}M_j^{o}] &= \left\{%
    \begin{array}{ll}
        \sum_{k=1}^N \Big(\rho_k \Psi_{v,k}^{({\bm 1}_i+{\bm 1}_j,0)}
        +\frac{\mathbb{E}[G_k]}{\mathbb{E}[C]}\Psi_{g,k}^{({\bm 1}_i +{\bm 1}_j,0)}
        +\frac{\mathbb{E}[S_k]}{\mathbb{E}[C]}\Psi_{s,k}^{({\bm 1}_i+{\bm 1}_j)}\Big)
        & \hbox{ if } i \neq j,\\
        2\sum_{k=1}^N \Big(\rho_k \Psi_{v,k}^{(2{\bm 1}_i,0)}
        +\frac{\mathbb{E}[G_k]}{\mathbb{E}[C]}\Psi_{g,k}^{(2{\bm 1}_i,0)}
        +\frac{\mathbb{E}[S_k]}{\mathbb{E}[C]}\Psi_{s,k}^{(2{\bm 1}_i)}\Big) & \hbox{ if } i=j,\\
    \end{array}%
    \right. \label{sec5-2-second-eqn1}\\
    \mathbb{E}[M_i^{oq}M_j^{oq}] &= \left\{%
    \begin{array}{ll}
        \mathbb{E}[M_i^{o}M_j^{o}]
        +\rho_i\Psi_{v,i}^{({\bm 1}_j,1)}+\rho_j\Psi_{v,j}^{({\bm 1}_i,1)}
        +\frac{\mathbb{E}[G_i]}{\mathbb{E}[C]}\Psi_{g,i}^{({\bm 1}_j,1)}
        +\frac{\mathbb{E}[G_j]}{\mathbb{E}[C]}\Psi_{g,j}^{({\bm 1}_i,1)} & \hbox{ if } i \neq j,\\
        \mathbb{E}[(M_i^{o})^2]
        +2\Big(\rho_i\Psi_{v,i}^{({\bm 1}_i,1)}
        +\frac{\mathbb{E}[G_i]}{\mathbb{E}[C]}\Psi_{g,i}^{({\bm 1}_i,1)}
        +\rho_i\Psi_{v,i}^{({\bm 0},2)}
        +\frac{\mathbb{E}[G_i]}{\mathbb{E}[C]}\Psi_{g,i}^{({\bm 0},2)}\Big) & \hbox{ if } i=j,\\
    \end{array}%
    \right. \label{sec5-2-second-eqn2}\\
    \mathbb{E}[M_iM_j] &= \left\{%
    \begin{array}{ll}
        \mathbb{E}[M_i^{oq}M_j^{oq}]+ \rho_i\Psi_{v,i}^{({\bm 1}_j,0)}+\rho_j\Psi_{v,j}^{({\bm 1}_i,0)}
        & \hbox{ if } i \neq j,\\
        \mathbb{E}[(M_i^{oq})^2]+ 2\rho_i\Psi_{v,i}^{({\bm 1}_i,0)}+2\rho_i\Psi_{v,i}^{({\bm 0},1)}
        +\rho_i  & \hbox{ if } i=j. \\
    \end{array}%
    \right. \label{sec5-2-second-eqn3}
\end{align}

\subsection{A numerical example}

In this subsection we present numerical results for the first and
second moments of the number of customers in each station. The
expression for the mean number of customers in each station is
given by (\ref{sec5-2-mean-eqn3}), together with
(\ref{sec5-2-mean-eqn1}) and (\ref{sec5-2-mean-eqn2}). By using
the formulas
(\ref{sec5-2-second-eqn1})-(\ref{sec5-2-second-eqn3}), we can
obtain an expression for the variance of the number of customers
in each station and an expression for the covariance of the
numbers of customers in two different stations. Note that these
moments are expressed in terms of $\Phi_i^{(\bm l,m)}$, refer to
Theorem \ref{sec5-new-thm1}. Therefore, these moments can be
obtained by using Theorems \ref{sec5-new-thm2} and
\ref{sec5-1-iteration}. In the following numerical example we
consider a single server polling model with five stations (i.e.,
$N=5$).

\noindent{\bf Example 1.~} We assume that the arrival rate of
type-$i$ customers is $\lambda_i=0.025$ for all $i$,
$i=1,\ldots,5$. The service times of type-$i$ customers are
exponentially distributed with means $\mathbb{E}[B_1]=1,
\mathbb{E}[B_2]=2, \mathbb{E}[B_3]=4, \mathbb{E}[B_4]=8$ and
$\mathbb{E}[B_5]=16$, respectively. Hence the total utilization of
the server is $\rho=\sum_{i=1}^5 \rho_i=0.775<1$. The switchover
times from station $i$ to station $i+1$ are deterministic with
$\mathbb{E}[S_i]=1$ for all $i$, $i=1,\ldots,5$. The retrial rate
of customers in the orbit of station $i$ is $\nu_i=1$ for all $i$,
$i=1,\ldots,5$. The glue periods at station $i$ are exponentially
distributed with parameters $\gamma_i$, $i=1,\ldots,5$. We assume
that $\gamma_i$ is the same for all $i$, i.e.,
$\mathbb{E}[G_i]=\mathbb{E}[G]$ for all $i$, $i=1,\ldots,5$.

 \begin{figure}
 \centering
 \subfigure[${0 \le \E[G] \le 10}$.]
 {\includegraphics[width=7cm]{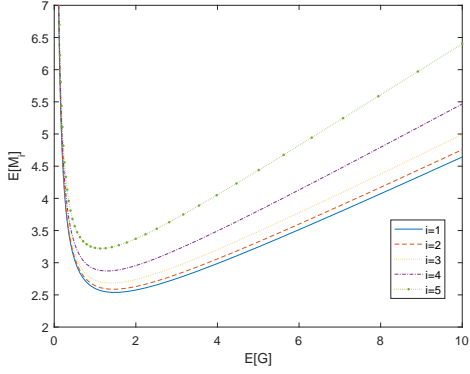}}
 \subfigure[${0 \le \E[G] \le 1000}$.]
 {\includegraphics[width=7cm]{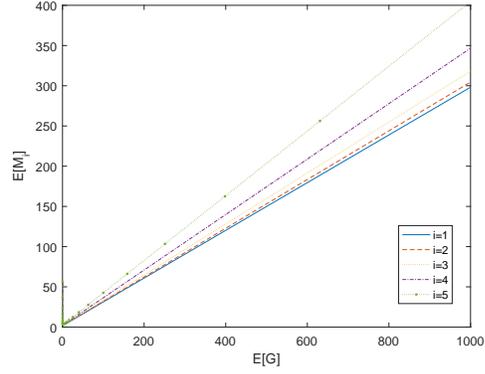} }
 \caption[Optional caption for list of figures]
 {The mean number of customers in station $i$, $\mathbb{E}[M_i]$, $i=1,\ldots,5$, varying $\E[G]$.} \label{fig1}
 \end{figure}

 \begin{figure}
 \centering
 \subfigure[${0 \le \E[G] \le 10}$.]
 {\includegraphics[width=7cm]{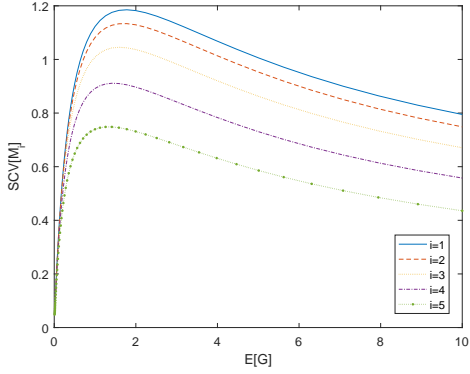}}
 \subfigure[${0 \le \E[G] \le 1000}$.]
 {\includegraphics[width=7cm]{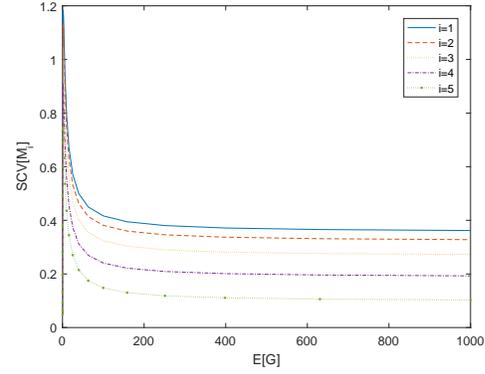} }
 \caption[Optional caption for list of figures]
 {The squared coefficient of variation for the number of customers in station $i$,
 $\mbox{SCV}[M_i]$, $i=1,\dots,5$, varying $\E[G]$.} \label{fig2}
 \end{figure}

 \begin{figure}
 \centering
 \subfigure[${0 \le \E[G] \le 10}$.]
 {\includegraphics[width=7cm]{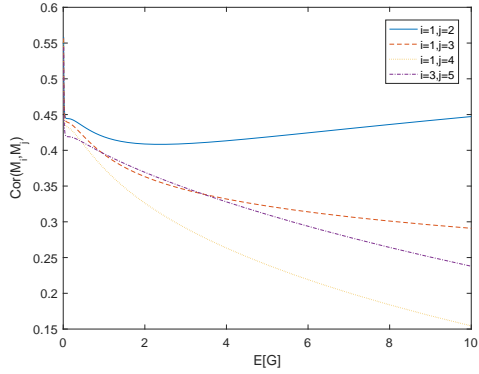}}
 \subfigure[${0 \le \E[G] \le 1000}$.]
 {\includegraphics[width=7cm]{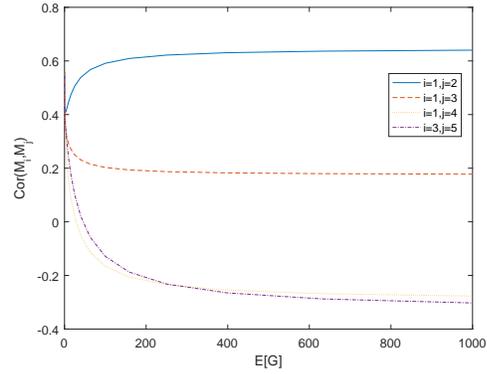} }
 \caption[Optional caption for list of figures]
 {The correlation coefficient of the numbers of customers in station $i$ and station $j$,
 $\mbox{Cor}(M_i,M_j)$, $(i,j)=(1,2)$, $(i,j)=(1,3)$, $(i,j)=(1,4)$ and $(i,j)=(3,5)$, varying $\E[G]$.} \label{fig3}
 \end{figure}

 In Figure \ref{fig1} we plot the mean number of customers in
station $i$, $\mathbb{E}[M_i]$, $i=1,\ldots,5$, varying the mean
glue period $\E[G]$. In Figure \ref{fig2} we plot the squared
coefficient of variation (SCV) for the number of customers in
station $i$, $\mbox{SCV}[M_i]$, $i=1,\ldots,5$, varying $\E[G]$.
In Figure \ref{fig3} we plot the correlation coefficient of the
numbers of customers in two different stations,
$\mbox{Cor}[M_1,M_2]$, $\mbox{Cor}[M_1,M_3]$,
$\mbox{Cor}[M_1,M_4]$ and $\mbox{Cor}[M_3,M_5]$, varying $\E[G]$.
In Figures \ref{fig1}(a), \ref{fig2}(a) and \ref{fig3}(a) we vary
$\E[G]$ from 0 to 10 in order to better reveal the behavior of the
system for small $\E[G]$. In Figures \ref{fig1}(b), \ref{fig2}(b)
and \ref{fig3}(b) we vary $\E[G]$ from 0 to 1000 in order to
examine the behavior of the system for large $\E [G]$.

We can draw the following conclusions from these plots:

\begin{itemize}
 \item For small glue period lengths, the chances for a customer to retry are very low, hence the station size is large.
 \item If the glue period is very large, the customers face a long delay before getting served.
 \item There exists an optimal glue length at which each station has a minimum mean station size.
\item The figures suggest that the following happens when the mean
glue period grows large:
\begin{enumerate}
 \item[(i)] The mean numbers of customers grow linearly in $\E [G]$.
\item[(ii)] The squared coefficient of variation tends to a limit
when $\E [G] \rightarrow \infty$. \item[(iii)] The correlation
coefficients between the numbers of customers in different
stations tend to some limit when $\E [G] \rightarrow \infty$.
\end{enumerate}
\end{itemize}

In \cite{Winands2011} the author considers classical polling
systems with a branching-type service discipline like exhaustive
or gated service, and without glue periods, for the case that
switchover times become large. It is readily seen that our polling
model starts to behave very similarly as such a polling model,
when the glue periods grow large; indeed, every type-$i$ customer
will now almost surely become glued during the first glue period
of station $i$ that it experiences during its stay in the system,
and hence will be served during the first visit period of station
$i$ after its arrival to the system - just as in an ordinary gated
polling system. However, we cannot immediately apply the
asymptotic results of \cite{Winands2011} where switchover times
become large, because it considers deterministic switchover times,
while the focus is on the waiting time distribution. In a future
paper we intend to study the asymptotic behaviour of polling
systems with large switchover times, thus also obtaining the
asymptotic behaviour of polling systems with large glue periods.
We shall, among others, derive asymptotic expressions for the
$k$th moment of the station size. Our preliminary findings are in
agreement with the limiting behaviours of the mean station size,
the squared coefficient of variation of the station sizes, and the
correlation coefficient of the station sizes. The mean station
size is asymptotically linear in the mean switchover times, and
the squared coefficient of variation and the correlation
coefficient of the station sizes converge as the mean switchover
times go to infinity (with the ratios of the swichover times being
constant), as displayed in the figures.

\section{The polling system with retrials and general glue periods}
\label{section_new}

In \cite{M_O_J_1} and in Section~\ref{section_old} of the current
paper we have presented the distribution and mean of the number of
customers at different time epochs for a gated polling model with
retrials and glue periods, where the glue periods are
deterministic and exponentially distributed, respectively. In this
section, we assume that glue periods have general distributions.
We first consider the distribution of the total workload in the
system and present a workload decomposition. Subsequently we use
this to obtain a pseudo conservation law, i.e., an exact
expression for a weighted sum of the mean waiting times. In its
turn, the pseudo conservation law is used to obtain an
approximation for the mean waiting times of all customer types. We
present numerical results that indicate that the approximation is
very accurate. Finally we use this approximation to optimize a
weighted sum of the mean waiting times, $\sum_{i=1}^N c_i
\E[W_i]$, where $c_i$, $i=1,\ldots,N$ are positive constants and
$\mathbb{E}[W_i]$ is the mean waiting time of a type-$i$ customer
until the start of its service, by choosing the glue period
lengths, given the total glue period in a cycle.

\subsection{Workload distribution and decomposition}

Define $V$ as the amount of work in the system in steady state.
Furthermore, let $\tilde B(s) = \sum_{i=1}^N \lambda_i (1- \tilde
B_i(s))$. The LST of the amount of work at an arbitrary time can
be written as
\begin{equation}
 \E[e^{-sV}] = \frac{1}{\E[C]} \sum_{i=1}^N\big(\E[S_i] \E [e^{-sV^{(S)}_{i}}] + \E[G_i] \E [e^{-sV^{(G)}_i}]+\rho_i \E[C] \E [e^{-sV_i^{(D)}}] \big),
 \label{work}
\end{equation}
where $V^{(S)}_{i}$, $V^{(G)}_{i}$ and $V^{(D)}_{i}$ are the
amount of work in the system during the switchover time from
station $i$, glue period of station $i$ and visit period of station $i$,
respectively.

Let $V_i^{(X)}$, $V_i^{(Y)}$ and $V_i^{(Z)}$ be the work in the
system at the \textit{start} of glue period of station $i$, visit period
of station $i$ and switchover period from station $i$, respectively. We know
that
\begin{align*}
\E [e^{-sV^{(S)}_{i}}] &= \E [e^{-sV_{i}^{(Z)}}]\frac{1-\tilde S_i(\tilde B(s))}{\E[S_i]\tilde B(s)}, \\
\E [e^{-sV^{(G)}_{i}}] &= \E [e^{-sV_{i}^{(X)}}]\frac{1-\tilde
G_i(\tilde B(s))}{\E[G_i]\tilde B(s)}.
\end{align*}
Therefore,
\begin{align}
& \sum_{i=1}^N \Big(\E [e^{-sV^{(S)}_{i}}] \E[S_{i}]+\E [e^{-sV^{(G)}_{i}}] \E[G_{i}]\Big) \nonumber \\
&= \sum_{i=1}^N\Bigg(\E [e^{-sV_i^{(Z)}}]\frac{1-\tilde S_i(\tilde
B(s))}{\tilde B(s)}+ \E [e^{-sV_{i}^{(X)}}]\frac{1-\tilde
G_i(\tilde B(s))}{\tilde B(s)}\Bigg)
 \nonumber \\
 &= \sum_{i=1}^N \Bigg(\frac{\E [e^{-sV_i^{(Z)}}]-\E [e^{-sV_i^{(Z)}}]\tilde S_i(\tilde B(s))+\E [e^{-sV_{i}^{(X)}}]-
 \E [e^{-sV_{i}^{(X)}}]\tilde G_i(\tilde B(s))}{\tilde B(s)}\Bigg) \nonumber \\
 &= \sum_{i=1}^N \Bigg(\frac{\E [e^{-sV_i^{(Z)}}]- \E [e^{-sV_{i+1}^{(X)}}] + \E [e^{-sV_i^{(X)}}] - \E [e^{-sV_{i}^{(Y)}}]}{\tilde B(s)}\Bigg)\nonumber \\
 &= \sum_{i=1}^N \Bigg(\frac{\E [e^{-sV_{i}^{(Z)}}]- \E [e^{-sV_{i}^{(Y)}}]}{\tilde B(s)}\Bigg).
 \label{nonvisitwork}
 \end{align}
Furthermore, using the last formula of the proof of Theorem 2 in
Boxma et al. \cite{BKK}, but with our notations, we have
 \begin{equation}
  \rho_i \E[C] \E [e^{-sV_i^{(D)}}] = \frac{\E [e^{-sV_{i}^{(Y)}}]- \E [e^{-sV_{i}^{(Z)}}]}{\tilde B(s)-s}.
 \label{visitwork}
 \end{equation}
Substituting \eqref{nonvisitwork} and \eqref{visitwork} in
\eqref{work}, we have
\begin{align}
 \E[e^{-sV}] = \frac{s}{\E[C](s-\tilde B(s))}\sum_{i=1}^N \Bigg(\frac{\E [e^{-sV_{i}^{(Z)}}]
 - \E [e^{-sV_{i}^{(Y)}}]}{\tilde B(s)}\Bigg). \label{work_final}
 \end{align}

Define the idle time as the time the server is not serving
customers (i.e., the sum of all the switchover and glue periods).
Let $V^{(Idle)}$ be the amount of work in the system at an
arbitrary moment in the idle time. We have, by
(\ref{nonvisitwork}) and (\ref{eqn-C}),
\begin{align}
\E [e^{-sV^{(Idle)}}] &=  \frac{1}{\E\big[\sum_{i=1}^N (S_i +  G_i)\big]}\sum_{i=1}^N
\Bigg(\E [e^{-sV^{(S)}_{i}}] \E[S_{i}]+\E [e^{-sV^{(G)}_{i}}] \E[G_{i}]\Bigg) \nonumber \\
 &= \frac{1}{\E\big[\sum_{i=1}^N (S_i +  G_i)\big]}\sum_{i=1}^N\Bigg(\frac{\E [e^{-sV_{i}^{(Z)}}]- \E [e^{-sV_{i}^{(Y)}}]}{\tilde B(s)}\Bigg) \nonumber\\
 &= \frac{1}{(1-\rho)\E[C]} \sum_{i=1}^N\Bigg(\frac{\E [e^{-sV_{i}^{(Z)}}]- \E [e^{-sV_{i}^{(Y)}}]}{\tilde B(s)}\Bigg).%\nonumber
\label{work_idle}
\end{align}

We know that the LST of the amount of work at steady state, $V_{M/G/1}$, in the standard $M/G/1$
queue where the arrival rate is $\sum_{i=1}^N \lambda_i$ and the LST of the service time distribution is
$\sum^N_{i=1} \frac{\lambda_i}{\sum^N_{j=1}\lambda_j} \tilde B_i(s)$, is given by
\begin{equation}
 \E[e^{-sV_{M/G/1}}] = \frac{(1-\rho)s}{s-\tilde B(s)}.
 \label{M/G/1}
\end{equation}
From Equations \eqref{work_final}, \eqref{work_idle} and
\eqref{M/G/1} we have
\begin{align*}
\E[e^{-sV}] = \E[e^{-sV_{M/G/1}}]\E[e^{-sV^{(Idle)}}].
\end{align*}

In Theorem 2.1 of \cite{Boxma89}, a workload decomposition property has been proved for
a large class of single-server multi-class
queueing systems with service interruptions (like switchover
periods or breakdowns). It amounts to the statement that, under
certain conditions, the steady-state workload is
in distribution equal to the sum of two independent quantities:
(i) the steady-state workload in the corresponding queueing model
without those interruptions, and (ii) the steady-state workload at
an arbitrary interruption epoch. The gated polling model with glue
periods and retrials of the present paper satisfies all the
assumptions of Theorem 2.1 of \cite{Boxma89}, and hence, in
agreement with what we have seen above, the workload decomposition
indeed holds.

\subsection{Pseudo conservation law}

By the workload decomposition, it is shown in \cite{Boxma89} that
  \begin{equation}
   \sum_{i=1}^N \rho_i \E [W_i] = \rho \frac{\sum_{i=1}^N \lambda_i \E [B_i^2]}{2(1-\rho)}
   + \rho \frac{\E\big[\big(\sum_{i=1}^N (S_i + G_i)\big)^2\big]}{2\E \big[\sum_{i=1}^N (S_i + G_i)\big]}
   + \frac{\E \big[\sum_{i=1}^N (S_i + G_i)\big]}{2(1-\rho)}\big(\rho^2-\sum_{i=1}^N \rho_i^2 \big) + \sum_{i=1}^N \E [F_i],
   \label{pcl_1}
  \end{equation}
where $F_i$ is the work left in station $i$ at the end of a visit
period of station $i$ (and hence at the start of a switchover from
station $i$). Other than $\E[F_i]$, Equation (\ref{pcl_1}) is
independent of the service discipline. Note that $ \E[F_i] =
\E[Z_i^{(i)}] \E[B_i]$. To find $\E[Z_i^{(i)}]$ we will derive a
relation between $\E[Z_i^{(i)}]$ and $\E[Y_i^{(iq)}]$.
$\E[Y_i^{(iq)}]$ consists of the following three parts:
\begin{itemize}
 \item[(i)] Mean number of type-$i$ customers who were already present at the end of the previous visit to station $i$
 and who are glued during the glue period just before the current visit to station $i$.
  \item[(ii)] Mean number of type-$i$ customers who have arrived during the time interval
  from the end of the previous visit to station $i$
  to the start of the glue period of station $i$ just before the visit to station $i$, and who are glued during that glue period.
  \item[(iii)] Mean number of type-$i$ customers who arrive during
  the glue period of station $i$ just before the visit to station $i$.
\end{itemize}
 Note that (i) equals $(1-\tilde G_i(\nu_i)) \E[Z_i^{(i)}]$ because the mean number of type-$i$ customers
 who were present at the end of the previous visit to station $i$ is $\E[Z_i^{(i)}]$, and
 the probability that a customer who was present at the end of the previous
 visit to station $i$ is glued during the glue period just before the current visit to station $i$, is
 $1-\tilde G_i(\nu_i)$.
 (ii) equals $(1-\tilde G_i(\nu_i)) \lambda_i \big((1-\rho_i) \E[C] - \E[G_i]\big)$.
 Here $\lambda_i \big((1-\rho_i) \E[C] - \E[G_i]\big)$ is the mean number of type-$i$ customers
 who have arrived during the time interval from the end of the previous visit to station $i$ to
 the start of the glue period of station $i$ just before the visit to station $i$.
 Finally, (iii) equals $\lambda_i \E[G_i]$.
 Therefore,
\begin{align}
\E [Y_i^{(iq)}]= & (1-\tilde G_i(\nu_i)) \E[ Z_i^{(i)}] + (1-\tilde G_i(\nu_i)) (1-\rho_i) \lambda_i \E [C]
+ \tilde G_i(\nu_i) \E [G_i] \lambda_i.
\label{JK}
\end{align}
Since $\rho_i=\frac{\E [Y_i^{(iq)}] \E[B_i]}{\E [C]}$, we have
$\E [Y_i^{(iq)}]=\lambda_i \E [C]$.
Hence, by (\ref{JK}), we get
 \begin{align*}
 \E[Z^{(i)}_i] &=\lambda_i \rho_i \E [C] + \frac{\lambda_i\tilde G_i(\nu_i)}{1-\tilde G_i(\nu_i)} \left(\E [C] - \E [G_i]\right).
 \end{align*}
Therefore, $\E[F_i]$ is given by
 \begin{equation}
 \E[F_i] =\rho_i^2 \E [C] + \frac{\rho_i \tilde G_i(\nu_i)}{1-\tilde G_i(\nu_i)} \left( \E [C] - \E [G_i]\right).
\label{moment1}
\end{equation}
The first term on the right-hand side equals the mean amount of
work for type-$i$ customers who arrived at station $i$ during a
visit period of station $i$. The second term is interpreted as
follows: Since $\lambda_i (\E [C]-\E [G_i])$ is the mean number of
type-$i$ customers who arrive during one cycle excluding the glue
period of station $i$ in that cycle, $(\tilde G_i(\nu))^k \rho_i
(\E [C]-\E [G_i])$ is the mean amount of work for type-$i$
customers who arrive during the $k$th previous cycle excluding the
glue period of station $i$ in that cycle, and who are present at
the end of the current visit period of station $i$. Hence, the
second term, which is $\sum^\infty_{k=1}(\tilde G_i(\nu))^k \rho_i
(\E [C]-\E [G_i])$, is the mean amount of work for type-$i$
customers who were present in the orbit of station $i$ at the
beginning of the visit period of station $i$.

From Equations \eqref{pcl_1} and \eqref{moment1} together with (\ref{eqn-C}), we obtain the following pseudo conservation law:
 \begin{align}
 \sum_{i=1}^N \rho_i \E W_i =& \rho \Bigg(\frac{\sum_{i=1}^N \lambda_i \E [B_i^2]}{2(1-\rho)}+
 \frac{\E \big[\big(\sum_{i=1}^N (S_i +  G_i)\big)^2\big]}{2 \E \big[\sum_{i=1}^N (S_i +
 G_i)\big]}\Bigg)
+ \Big(\rho^2 + \sum_{i=1}^N \rho_i^2 \Big)\frac{\E \big[\sum_{i=1}^N (S_i +  G_i)\big]}{2(1-\rho)} \nonumber \\
&+ \sum_{i=1}^N \frac{\rho_i \tilde G_i(\nu_i)}{1-\tilde
G_i(\nu_i)}\Bigg(\frac{\E \big[\sum_{j=1}^N (S_j +
G_j)\big]}{1-\rho} - \E [G_i] \Bigg). \label{pcl_2}
\end{align}

\subsection{Approximation of the mean waiting times}

We now use the pseudo conservation law to find an approximation
for the mean waiting times of all customer types. Below we briefly
sketch the idea behind the approximation. Everitt \cite{Everitt}
has developed a method to approximate the mean waiting times in an
ordinary gated polling system (without retrials and glue periods).
The idea in this approximation is that an arriving customer first
has to wait for the residual cycle time, until the server begins a
new visit to its station. Subsequently, it has to wait for the
service times of all customers of the same type, who arrived
before it, in the elapsed cycle time. This leads to $\E [W_j] =
(1+ \rho_j) \E [R_{c_j}]$, where $\E [R_{c_j}]$ is the mean of the
residual time of a cycle starting with a visit to station $j$,
which is the same as the mean of the elapsed time of a cycle
starting with a visit to station $j$. Next, Everitt assumed that
for all $j$, this mean residual cycle time is independent of $j$,
i.e., $\E [R_{c_j}] \approx \E[R_{c}]$, leading to the
approximation $\E [W_j] \approx (1+ \rho_j) \E [R_{c}]$ for the
model without retrials and glue periods.

In this paper, we introduce a similar type of approximation, including
one extra term, for the mean waiting times in the model with retrials and glue periods:
\begin{equation}
    \E [W_j] \approx (1+\rho_j) \E[R_c] + \frac{\tilde G_j (\nu_j)}{1-\tilde G_j (\nu_j)} \left(\E [C] - \E [G_j]\right).
    \label{approx_time_}
\end{equation}
In the appendix we provide a detailed derivation of
(\ref{approx_time_}). The first term on the right-hand side of
(\ref{approx_time_}) is the same as the term in \cite{Everitt}.
The second term on the right-hand side is added because not every
customer who arrives in a particular cycle receives service in
that cycle. The type-$j$ customers arriving during any period
other than the glue period of station $j$ receive service in the
following visit period with probability $1-\tilde{G_j}(\nu_j)$.
Furthermore, the type-$j$ customers arriving during any period
other than a glue period of station $j$, have to wait for a
geometric number (with parameter $1-\tilde{G_j}(\nu_j)$) of cycles
before receiving service. Since a type-$j$ customer arrives during
a period other than a glue period of station $j$ with probability
$\frac{\E[C]-\E[G_j]}{\E[C]}$, the mean number of cycles until an
arbitrary type-$j$ customer receives its service, is
$\frac{\E[C]-\E[G_j]}{\E[C]}\times \frac{\tilde G_j
(\nu_j)}{1-\tilde G_j (\nu_j)}$. The second term is obtained by
multiplying this mean number of cycles and the mean cycle time.

It should be noted that, in reality, the mean residual cycle times
for station $i$ and station $j$ ($j\neq i$) are not equal. A key
element of our approximation is to assume that they {\em are}
equal. We can now use the pseudo conservation law to determine the
one unknown term $\E [R_c]$: By substituting \eqref{approx_time_}
in \eqref{pcl_2} and using (\ref{eqn-C}), we get
%\begin{align*}
%   & \sum_{i=1}^{N} \rho_i(1+\rho_i) \E[R_c]+ \sum_{i=1}^N \frac{\rho_i \tilde G_i(\nu_i)}{1-\tilde G_i(\nu_i)}\left(\E[C] - \E [G_i] \right) \\
%   &\approx
%    \rho \Bigg(\frac{\sum_{i=1}^N \lambda_i \E[B_i^2]}{2(1-\rho)}
%    + \frac{\E \big[\big(\sum_{i=1}^N (S_i+G_i)\big)^2\big]}{2 \E \big[\sum_{i=1}^N (S_i +  G_i)\big]}\Bigg)
%    + \Big(\rho^2 + \sum_{i=1}^N \rho_i^2 \Big)\frac{\E \big[\sum_{i=1}^N (S_i +  G_i)\big]}{2(1-\rho)} \nonumber\\
%    & \quad + \sum_{i=1}^N \frac{\rho_i \tilde G_i(\nu_i)}{1-\tilde G_i(\nu_i)}
%    \Bigg(\frac{\E\big[\sum_{j=1}^N (S_j+G_j)\big]}{1-\rho}-\E [G_i]\Bigg).
%\end{align*}
\begin{align}
    \E[R_c] &\approx
    \frac{\rho}{\rho + \sum_{i=1}^{N} \rho_i^2} \Bigg(\frac{\sum_{i=1}^N \lambda_i \E [B_i^2]}{2(1-\rho)}
    + \frac{\E \big[\big(\sum_{i=1}^N (S_i+G_i)\big)^2\big]}{2 \E \big[\sum_{i=1}^N (S_i +  G_i)\big]}
    + \rho \frac{\E\big[\sum_{i=1}^N (S_i +  G_i)\big]}{2(1-\rho)}
    \Bigg) \nonumber \\
    &\quad+\frac{\sum_{i=1}^N \rho_i^2}{\rho+\sum_{i=1}^{N} \rho_i^2} \Bigg(\frac{\E \big[\sum_{i=1}^N(S_i+G_i)\big]}{2(1-\rho)}\Bigg).
    \label{residual_cyc}
\end{align}
Substitution of \eqref{residual_cyc} into \eqref{approx_time_}
yields the approximation for the mean waiting times of all
customer types:
\begin{align}
    \E [W_j] &\approx \frac{1+\rho_j}{\rho + \sum_{i=1}^{N}
    \rho_i^2}\Bigg\{
    \rho \Bigg(\frac{\sum_{i=1}^N \lambda_i \E [B_i^2]}{2(1-\rho)}
    + \frac{\E \big[\big(\sum_{i=1}^N (S_i +  G_i)\big)^2\big]}{2 \E \big[\sum_{i=1}^N (S_i +  G_i)\big]}
    + \rho \frac{\E \big[\sum_{i=1}^N (S_i +  G_i)\big]}{2(1-\rho)}\Bigg) \nonumber\\
    &\quad+ \sum_{i=1}^N \rho_i^2\Bigg(\frac{\E \big[\sum_{i=1}^N (S_i +  G_i)\big]}{2(1-\rho)}\Bigg) \Bigg\}
     + \frac{\tilde G_j (\nu_j)}{1-\tilde G_j (\nu_j)} \left(\frac{\E \big[\sum_{i=1}^N (S_i+ G_i)\big]}{1-\rho}-\E [G_j]\right).
    \label{approx_time_final}
\end{align}

We now consider various examples to compare the above
approximation results with the exact analysis from \cite{M_O_J_1}
for deterministic glue periods and from Section \ref{section_old}
of the present paper for exponentially distributed glue periods.
Further, we will compare the results of this approximation with
simulation results for the case that the glue periods follow a
gamma distribution.

\subsection*{Deterministic glue periods}
In the numerical example of Table \ref{tab:1} we consider a
two-station polling system. The switchover times and service times
are exponentially distributed. We keep the parameters of station
$1$ fixed, $\lambda_1= 1, \E[B_1] =0.45, \E[S_1]=1, G_1=0.5, \nu_1
= 1$, and vary the parameters of station $2$.

\subsection*{Exponential glue periods}
In the numerical example of Table~\ref{tab:2} we consider a
three-station polling system. The switchover times are
deterministic, and the service times are exponentially
distributed. We keep the parameters of station $1$ fixed,
$\lambda_1= 1, \E[B_1] =0.45, \E[G_1]=0.5$. Further, the
switchover times and exponential retrial rates of all three
stations are fixed, $S_1 =S_2 =S_3 =1$ and $\nu_1=\nu_2=\nu_3 =
1$.

\begin{table}[t]
\begin{center}
\scriptsize
  \begin{tabular}{| c | c | c | c | c | c | c |}
  \hline
  $\lambda_2$ & $\E[B_2]$ & $\E[S_2]$ & $G_2$ & $\nu_2$ & Exact $(\E[W_1], \E[W_2])$ & Approx $(\E[W_1], \E[W_2])$ \\ \hline
  1&0.45&1&0.5&1& (71.61, 71.61) &(71.61, 71.61) \\ \hline
  0.5&0.45&1& 0.5&1& (21.44, 20.34) & (21.49, 20.24)\\ \hline
  0.5&0.2&1& 0.5&1&(15.18, 13.96) & (15.21, 13.83)\\ \hline
  0.5&0.2&2& 0.5&1& (20.52, 18.82) & (20.55, 18.71)\\ \hline
  0.5&0.2&2& 1&1& (23.01, 11.48) & (22.99, 11.67)\\ \hline
  0.5&0.2&2& 1&0.5& (22.97, 20.31) & (22.99, 20.20)\\ \hline
  \end{tabular}
 \caption{Comparison of exact and approximate mean waiting times for a polling system with deterministic glue periods.}
          \label{tab:1}
 \end{center}
 \end{table}

\begin{table}[t]
\begin{center}
\scriptsize
\begin{tabular}{| c | c | c | c | c | c | c | c |}
     \hline
  $\lambda_2$ & $\E[B_2]$ & $\E[G_2]$ & $\lambda_3 $ &  $\E[B_3]$ & $\E[G_3]$ & Exact $(\E[W_1], \E[W_2],
  \E[W_3])$ & Approx $(\E[W_1], \E[W_2], \E[W_3])$  \\ \hline
  1&0.3&0.5&1&0.3&0.5 & (121.0, 121.0, 121.0)  & (121.0, 121.0, 121.0)\\ \hline
  1&0.3&0.5&0.5&0.3&0.5 & (47.59, 47.58, 46.74) & (47.71, 47.71, 46.24)\\ \hline
  1&0.3&0.5&0.5&0.1&0.5 & (33.65, 33.64, 32.54) & (33.69, 33.69, 31.97)\\ \hline
  2&0.3&0.5&0.5&0.1&0.5 & (246.8, 246.6,  242.3) & (242.4, 257.1, 230.2)\\ \hline
  2&0.15&0.5&0.5&0.1&0.5 & (33.52, 33.51, 32.42) & (33.56, 33.56, 31.86)\\ \hline
  2&0.15&2&0.5&0.1&0.5 & (44.88, 19.71, 43.64) & (45.22, 19.50, 42.92)\\ \hline
  2&0.15&2&0.5&0.1&1 & (48.66, 21.42, 28.75) & (49.03, 21.17, 27.98)\\ \hline
  \end{tabular}
 \caption{Comparison of exact and approximate mean waiting times for a polling system with exponentially distributed glue periods.}
   \label{tab:2}
   \end{center}
  \end{table}

\subsection*{Gamma distributed glue periods}
In the above two examples we can get the exact mean waiting times
using the method in \cite{M_O_J_1} and Section \ref{section_old}
of this paper. In the numerical examples of Table \ref{tab:4} we
compare the approximate mean waiting times with simulation
results, for a polling system where the lengths of glue periods
are gamma distributed.
% therefore Table A (on page 27) continues/follows onto Table B (on page 28).

We consider a five-station polling system in which the glue
periods, switchover times and service times are all gamma
distributed. We simulate such a system to find the mean waiting
times. We also give a 95\% confidence interval for the mean
waiting times obtained using simulations. We have generated one
million cycles, splitting this into ten periods of $10^5$ cycles,
and using the results of these ten periods to obtain confidence
intervals. Then we compare the simulation results with the results
obtained using the approximation formula. Here, $k$ and $\theta$
are, respectively, the shape and the scale parameters of the gamma
distribution with probability density function $\frac{1}{\Gamma(k)
\theta^k} x^{k-1} e^{- \frac{x}{\theta}}$.

\subfigcapskip 3 mm
\begin{table}[t]
\begin{center}
\subtable[The parameter values for the four cases.]{
\scriptsize
  \begin{tabular}{| c | c | c| c| c| c|c|c|c|}
  \hline
  Parameters & (i) & (ii) & (iii) & (iv) \\ \hline
   ($\lambda_1$, $\lambda_2$, $\lambda_3$, $\lambda_4$, $\lambda_5$) &(0.1, 0.1, 0.1, 0.1, 0.1) &(0.1, 0.1, 0.1, 0.1,
     0.2)
     &(0.1, 0.1, 0.1, 0.3, 0.2) &(0.1, 0.2, 0.1, 0.2, 0.2)\\ \hline
     ($\nu_1$, $\nu_2$, $\nu_3$, $\nu_4$, $\nu_5$) & (2.0, 2.0, 2.0, 2.0, 2.0) &(2.0, 2.0, 2.0, 2.0, 3.0) &(2.0, 2.0, 2.0, 1.0,
     3.0)
     &(2.0, 5.0, 2.0, 4.0, 3.0) \\ \hline
     ($k_{B_1}$, $k_{B_2}$, $k_{B_3}$, $k_{B_4}$, $k_{B_5}$) &(1.0, 1.0, 1.0, 1.0, 1.0) &(1.0, 1.0, 1.0, 1.0,
     0.8)
     &(1.0, 1.0, 1.0, 0.5, 0.8) &(1.0, 0.5, 1.5, 0.5, 0.8)\\ \hline
     ($\theta_{B_1}$, $\theta_{B_2}$, $\theta_{B_3}$, $\theta_{B_4}$, $\theta_{B_5}$) &(1.5, 1.5, 1.5, 1.5,
     1.5)
     &(1.5, 1.5, 1.5, 1.5, 1.5) &(1.5, 1.5, 1.5, 1.5, 1.5) &(1.5, 1.5, 1.5, 1.5, 1.5) \\ \hline
     ($k_{S_1}$, $k_{S_2}$, $k_{S_3}$, $k_{S_4}$, $k_{S_5}$) &(2.0, 2.0, 2.0, 2.0, 2.0) &(2.0, 2.0, 2.0, 2.0,
     1.0)
     &(2.0, 2.0, 2.0, 3.0, 1.0) &(2.0, 5.0, 2.0, 3.0, 1.0) \\ \hline
     ($\theta_{S_1}$, $\theta_{S_2}$, $\theta_{S_3}$, $\theta_{S_4}$, $\theta_{S_5}$) &(1.0, 1.0, 1.0, 1.0,
     1.0)
     &(1.0, 1.0, 1.0, 1.0, 1.0) &(1.0, 1.0, 1.0, 1.0, 1.0) &(1.0, 1.0, 1.0, 1.0, 1.0) \\ \hline
     ($k_{G_1}$, $k_{G_2}$, $k_{G_3}$, $k_{G_4}$, $k_{G_5}$) &(1.0, 1.0, 1.0, 1.0, 1.0) &(1.0, 1.0, 1.0, 1.0,
     2.0)
     &(1.0, 1.0, 1.0, 0.5, 2.0) &(1.0, 3.0, 1.0, 0.5, 2.0)\\ \hline
     ($\theta_{G_1}$, $\theta_{G_2}$, $\theta_{G_3}$, $\theta_{G_4}$, $\theta_{G_5}$) &(1.0, 1.0, 1.0, 1.0,
     1.0)
     &(1.0, 1.0, 1.0, 1.0, 1.0) &(1.0, 1.0, 1.0, 1.0, 1.0) & (1.0, 1.0, 1.0, 1.0, 1.0) \\ \hline
  \end{tabular}
  }
 \subtable[Simulation results.]{
\scriptsize
  \begin{tabular}{| c | c | c|c|}
  \hline
  Cases & ($\E [W_1], \E [W_2], \E[W_3], \E[W_4],
  \E[W_5]$) & 95\% lower confidence bound & 95\% upper confidence
  bound \\ \hline
  (i) & (68.94, 68.92, 68.87, 68.91, 68.84) & (68.65, 68.64, 68.60, 68.58, 68.54) & (69.23, 69.21, 69.14, 69.24, 69.15) \\ \hline
  (ii) & (108.59, 108.47, 108.40, 108.28, 72.43)  & (107.86, 107.66, 107.58, 107.50, 71.99) & (109.33, 109.27, 109.23, 109.06, 72.87) \\ \hline
  (iii) & (217.54, 218.44, 219.51, 548.58, 144.61) & (216.04, 217.11, 218.01, 544.87, 143.68) & (219.04, 219.77, 221.01, 552.28, 145.54) \\ \hline
  (iv) & (276.39, 158.43, 283.36, 343.52, 183.40)  & (275.30, 157.84, 282.15, 342.07, 182.71) & (277.47, 159.01, 284.57, 344.97, 184.09)\\
  \hline
\end{tabular}
 }
 \subtable[Approximation results.]{
 \scriptsize
  \begin{tabular}{| c | c |}
  \hline
  Cases & Approx ($\E [W_1], \E  [W_2],\E[W_3],\E[W_4],\E[W_5]$) \\
  \hline
 (i) & (69.00, 69.00, 69.00, 69.00, 69.00) \\ \hline
 (ii) & (108.25, 108.25, 108.25, 108.25, 72.84) \\ \hline
 (iii) & (210.22, 210.22, 210.22, 566.37, 140.93) \\ \hline
  (iv) & (274.10, 155.12, 284.14, 348.72, 182.01) \\ \hline
  \end{tabular}
 }
 \end{center}
 \vspace{-0.3cm}
 \caption{Comparison of mean waiting times from simulations and approximate mean waiting times
 for a polling system with gamma distributed glue periods.}
 \label{tab:4}
 \end{table}

The values of the parameters are listed in Table \ref{tab:4}(a).
Table \ref{tab:4}(b) shows the mean waiting times by simulation,
along with 95\% lower and upper confidence bounds. Table
\ref{tab:4}(c) shows the approximate mean waiting times. We can
draw the following conclusions about the mean waiting time
approximation.
\begin{itemize}
 \item The mean waiting time approximation is very accurate. In only two cases
 (the fourth case in Table \ref{tab:2} and the case (iii) in Table \ref{tab:4})
the error is in the order of $5\%$; in all other cases, we find errors which typically are less than $2 \%$.
 \item The mean waiting time approximation at one station is independent of the change in retrial rates of other stations,
 which is not true in reality.
 \item The mean waiting time approximations for two totally symmetric stations are the same,
 independent of their order in the system; but this is also not quite true in reality.
 \end{itemize}

 \subsection{Optimal choice of the glue period distributions} \label{section_optimise}
 In this subsection we discuss an optimization problem for the choice of the distributions
 of the glue periods, $G_i$, $i=1,\ldots,N$, to minimize the weighted sum of the mean waiting times
 $\sum^N_{i=1} c_i \E[W_i]$, subject to the constraint $\sum^N_{i=1} \E[G_i]=L$,
 where $c_i$, $i=1,\ldots,N$, and $L$ are positive constants.
 Because we do not have an explicit formula for the mean waiting time, it is difficult to solve exactly
 the constrained minimization problem.
 Instead of finding the exact solution of the constrained minimization problem,
 we will find the optimal choice of the distributions of $G_i$, $i=1,\ldots,N$,
 to
 \begin{align*}
 & \mbox{minimize} \qquad \sum^N_{i=1} c_i U_i \\
 & \mbox{subject to} \quad ~~ \sum^N_{i=1} \E[G_i] =L,
 \end{align*}
 where $U_i$ is the approximation of $\E[W_i]$ given by the right-hand side of (\ref{approx_time_final}).
 Note that under the constraint $\sum^N_{i=1} \E[G_i] =L$, the objective function
 of the minimization problem becomes
 \begin{align}
 \sum^N_{i=1} c_i U_i =& \sum^N_{i=1} \frac{c_i(1+\rho_i)}{\rho+\sum_{j=1}^N
    \rho_j^2}\Bigg[
 \rho \Bigg(\frac{\sum_{j=1}^N \lambda_j \E [B_j^2]}{2(1-\rho)}+ \frac{\E[(\sum_{j=1}^N S_j)^2]+2L
 \sum^N_{j=1}\E[S_j]
    + \E[(\sum^N_{j=1}G_j)^2]}{2 \E [\sum_{j=1}^N S_j +  L]} \Bigg) \nonumber \\
 & + \frac{\E [\sum_{j=1}^N S_j +  L]}{2(1-\rho)}\Big(\rho^2 +
 \sum_{j=1}^N \rho_j^2 \Big) \Bigg] + \sum^N_{i=1}\frac{c_i\E[e^{-\nu_i G_i}]}{1-\E[e^{-\nu_i
        G_i}]} \left(\frac{\E\big[\sum_{j=1}^N S_j + L\big]}{1-\rho}
 - \E[G_i]\right). \label{Opt-1}
 \end{align}
 By Jensen's inequality, it can be shown that if the
 nondeterministic glue period distributions with  means $g_i$, are changed to
 the degenerate (deterministic) ones with the same means $g_i$, $i=1,\dots,N$, then the
 right-hand side of (\ref{Opt-1}) becomes strictly smaller. Therefore,
 the above optimization problem becomes as follows:
 \begin{align}
 & \mbox{minimize}  \qquad  U(g_1,\ldots,g_N) \nonumber\\
 & \mbox{subject to} \nonumber\\
 & \qquad \qquad \qquad g_i>0, \quad i=1,\ldots, N, \label{con1}
 \\
 & \qquad \qquad \qquad \sum^N_{i=1}g_i =L, \label{con2}
 \end{align}
 where
 \begin{align}
 & U(g_1,\ldots,g_N) \nonumber\\
 &= \sum^N_{i=1} \frac{c_i(1+\rho_i)}{\rho + \sum_{j=1}^{N}
    \rho_j^2}\Bigg[
 \rho \Bigg(\frac{\sum_{j=1}^N \lambda_j \E [B_j^2]}{2(1-\rho)}
 + \frac{\E[(\sum_{j=1}^N S_j)^2]+2L \sum^N_{j=1}\E[S_j]
    +L^2]}{2 \E [\sum_{j=1}^N S_j +  L]}\Bigg)
 \nonumber\\
 & \quad + \frac{\sum_{j=1}^N \E[S_j] +  L}{2(1-\rho)}\Big(\rho^2 + \sum_{j=1}^N \rho_j^2 \Big) \Bigg]
 + \sum^N_{i=1} c_i \Big(-1 + \frac{1}{1-e^{-\nu_i g_i}}\Big) \Bigg(\frac{\sum_{j=1}^N \E[S_j]+L}{1-\rho} - g_i\Bigg).
\label{opt_res}
 \end{align}
 Since $U(g_1,\ldots,g_N)$ is continuous on ${\cal D} \equiv \{(g_1,\ldots,g_N): g_1>0,\ldots,g_N>0, g_1+\cdots+g_N=L\}$
 and $U(g_1,\ldots,g_N)\to \infty$ as $\min \{g_1,\ldots,g_N\}\to 0+$,  $U(g_1,\ldots,g_N)$ takes
 a minimum at a point in ${\cal D}$.
 At a minimum point $(g_1,\ldots,g_N)$, there exists a Lagrange multiplier $\kappa$ satisfying
 \begin{align}
 f_i(g_i)&=\kappa, \quad i=1,\ldots,N, \label{Larg}
 \end{align}
 where
 \begin{align*}
 f_i(g_i) &\equiv
 c_i-\frac{c_i}{1-e^{-\nu_i g_i}}
 -\frac{c_i\nu_ie^{-\nu_i g_i}}{(1-e^{-\nu_i g_i})^2}\Bigg(\frac{\sum_{j=1}^N \E[S_j] +
    L}{1-\rho}-g_i\Bigg), \quad i=1,\ldots,N.
 \end{align*}
 For each $i=1,\ldots,N$, the function $f_i:(0,L) \to (-\infty,f_i(L))$ is bijective,
 continuous and strictly increasing. Therefore, it has the inverse function
 $h_i:(-\infty,f_i(L)) \to (0,L)$, which is also continuous and
 strictly increasing. Therefore, Equation (\ref{Larg}) and the
 constraints (\ref{con1})  and (\ref{con2}) can be written as
 \begin{align}
 & \sum^N_{j=1} h_j(\kappa)= L, \quad -\infty< \kappa < \min\{f_1(L),\ldots,f_N(L)\}, \label{sol1} \\
 & g_i = h_i(\kappa), \quad i=1,\ldots, N. \label{sol2}
 \end{align}
 Since $\lim_{\kappa\to-\infty}\sum^N_{j=1} h_j(\kappa)=0$,
 $\lim_{\kappa \to (\min\{f_1(L),\ldots,f_N(L)\})-} \sum^N_{j=1} h_j(\kappa) > L$
 and $\sum^N_{j=1} h_j(\kappa)$ is strictly increasing in $\kappa$,
 (\ref{sol1}) has a unique solution, say $\kappa^*$.
 Therefore, from (\ref{sol2}), the optimal solution $(g_1^*,\ldots,g_N^*)$ is given
 by
 \begin{align*}
 g_i^* = h_i(\kappa^*), \quad i=1,\ldots, N.
 \end{align*}

 We will now consider a few numerical examples to look at the dependency of different system characteristics
 and the respective optimal glue periods.
 In \cite{M_O_J_Ton} a similar system was studied with a focus on optical switches,
 where the revenue of the system depended on distributing glue periods optimally to each station.
 In these examples we will look at the problem of minimizing $\sum_{i=1}^N c_i U_i$, that is the weighted waiting cost
 of the system given that the sum of expected values of glue periods is fixed. Since the optimization
 problem showed that the system performs best when the glue periods are deterministic,
 we will only consider models with deterministic glue periods.

We consider a three-station model and in each case vary one
parameter to study how the system performs under certain changes.
In all the cases the sum of the lengths of deterministic glue
periods is fixed, $L= 3$, and the service times and the switchover
times are exponentially distributed. The switchover times are
symmetric and fixed for all three stations, i.e. $\E[S_i] = 2$ for
all $i = 1,2,3$.

\begin{itemize}
\item[(i)] Case 1: In this case we keep all system parameters
symmetric except the arrival rate $\lambda_i$ of each station. Let
$\nu_i=1$, $\E[B_i]=1$ and $c_i=1$ for all $i=1,2,3$. In Table
\ref{tab:5} we show the optimal values of $g_1, g_2, g_3$ and
$\sum_{i=1}^N c_i U_i$ for different values of $\lambda_i$.

\item[(ii)] Case 2: In this case we keep all system parameters
symmetric except the mean service time $\E[B_i]$ of each station.
Let $\lambda_i= 1$, $\nu_i=1$, and $c_i=1$ for all $i=1,2,3$. In
Table \ref{tab:6} we show the optimal values of $g_1, g_2, g_3$
and $\sum_{i=1}^N c_i U_i$ for different values of $\E[B_i]$.

\item[(iii)] Case 3: In this case we keep all system parameters
symmetric except the retrial rate $\nu_i$ of each station. Let
$\lambda_i=0.25$, $\E[B_i]=1$ and $c_i=\rho_i$ for all $i=1,2,3$.
In Table \ref{tab:7} we show the optimal values of $g_1, g_2, g_3$
and $\sum_{i=1}^N c_i U_i$ for different values of $\nu_i$. Note
that in this case $\sum_{i=1}^N c_i U_i = \sum_{i=1}^N \rho_i
\E[W_i]$.

\item[(iv)] Case 4: In this case we keep all system parameters
symmetric except the weight $c_i$ of each station. Let
$\lambda_i=0.25$, $\E[B_i]=1$ and $\nu_i=1$ for all $i=1,2,3$. In
Table \ref{tab:8} we show the optimal values of $g_1, g_2, g_3$
and $\sum_{i=1}^N c_i U_i$ for different values of $c_i$.
\end{itemize}

 We can draw the following conclusions about the optimal allocation of glue periods using the above method.
\begin{itemize}
 \item The allocation doesn't depend on the arrival rate or mean service time of a station.
 This is due to the following observation:
 The first term on the right-hand side of Equation \eqref{opt_res} is independent of $g_i$
 and the second term is independent of
 arrival rates and mean service times. This might not be the case in exact analysis.
 \item The higher the retrial rate, the shorter the length of the glue period assigned to the station.
 \item The higher the weight allocated to a station, the bigger the length of the glue period assigned to the station.
 This helps us in scenarios when a waiting cost is associated with stations.
\end{itemize}

\begin{table}[t]
\begin{center}
\begin{tabular}{| c |  c | c | c | c  | c  |c|}
     \hline
  $\lambda_1$ &  $\lambda_2$    & $\lambda_3$ & $g_1 $ &  $g_2$   & $g_3$ & $\sum_{i=1}^N c_i U_i$   \\ \hline
  0.3&0.3&0.3& 1&1&1& 359.898\\ \hline
  0.3&0.2&0.2&1&1&1&115.063  \\ \hline
  0.3&0.2&0.1&1&1&1&84.362 \\ \hline
  \end{tabular}
  \caption{Optimal length of glue periods for different arrival rates.}
   \label{tab:5}
 \end{center}
\end{table}

\begin{table}[t!]
\begin{center}
\begin{tabular}{| c |  c | c | c | c  | c  |c|}
     \hline
  $\E[B_1]$ &  $\E[B_2]$  & $\E[B_3]$ & $g_1 $ &  $g_2$   & $g_3$ & $\sum_{i=1}^N c_i U_i$   \\ \hline
  0.3&0.3&0.3& 1&1&1& 422.888\\ \hline
  0.3&0.2&0.2&1&1&1&137.887  \\ \hline
  0.3&0.2&0.1&1&1&1&101.876 \\ \hline
  \end{tabular}
  \caption{Optimal length of glue periods for different mean service times.}
   \label{tab:6}
\end{center}
\end{table}

\begin{table}[t!]
\begin{center}
\begin{tabular}{| c |  c | c | c | c  | c  |c|}
     \hline
  $\nu_1$ &  $\nu_2$    & $\nu_3$ & $g_1 $ &  $g_2$   & $g_3$ & $\sum_{i=1}^N c_i U_i$   \\ \hline
  3&3&3& 1.0000&1.0000&1.0000& 84.001\\ \hline
  3&2&2&0.8340&1.0830&1.0830&90.680  \\ \hline
  3&2&1&0.7134&0.9157&1.3710&101.679 \\ \hline
  \end{tabular}
  \caption{Optimal length of glue periods for different retrial rates.}
   \label{tab:7}
   \end{center}
\end{table}

\begin{table}[t!]
\begin{center}
 \begin{tabular}{| c |  c | c | c | c  | c  |c|}
     \hline
  $c_1$ &  $c_2$ & $c_3$ & $g_1 $ &  $g_2$   & $g_3$ & $\sum_{i=1}^N c_i U_i$   \\ \hline
  3&3&3& 1.0000&1.0000&1.0000& 418.823\\ \hline
  3&2&2&1.1268&0.9366&0.9366&323.736  \\ \hline
  3&2&1&1.2311&1.0263&0.7426&271.086 \\ \hline
  \end{tabular}
  \caption{Optimal length of glue periods for different weights.}
   \label{tab:8}
\end{center}
\end{table}

 \section{Suggestions for further research} \label{section_concl}

 In this paper we have studied a gated polling model with the special features of retrials
 and glue, or reservation periods. For the case of exponentially distributed
glue periods, we have presented an algorithm to obtain the moments
of the number of customers in each station. We would like to point
out that phase-type glue periods can in principle be handled by
the same method.

For generally distributed glue periods, we have obtained an expression for the steady-state
distribution of the total workload in the system,
and we have used it to derive a pseudo conservation law for a weighted sum of the mean waiting times,
which in turn led us to an accurate approximation of the individual mean waiting times.
A topic for further research is to analyze the exact waiting time distribution,
for exponentially distributed glue periods and for constant glue periods.

The introduction of the concept of glue period was motivated by
the wish to obtain insight into the performance of certain
switches in optical communication systems. We have considered the
optimal choice of the glue period lengths, under the constraint
that the total glue period length per cycle is fixed. A topic for
further study is the unconstrained counterpart to this
optimization problem; a complication one then faces is that the
objective function for the optimization can be nonconvex. In fact,
it is possible that the Hessian of the objective function is not
positive semi-definite even for the two-station system. However it
still seems to be intuitively natural that there will exist a
unique solution for the optimization problem.

Not restricting ourselves to optical communications, one can also
interpret a glue period as a reservation period - a window of
opportunity for claiming service at the next visit of the server
to a station. It would be interesting to study reservation periods
in more detail, and in particular to consider the problem of
choosing reservation periods in such a way that some objective
function is optimized.

\appendix
\section*{Appendix: Approximation of the mean waiting times}

Below we outline a method to approximate the mean waiting times of
all customer types. The arrival of a type-$i$ customer occurs
either during a glue period of station $i$ or during any other
period. At the start of the visit period the customers which will
be served in the current visit period are fixed. The mean length
of the visit period is now the same irrespective of the order in
which these customers are served. Without loss of generality, we
will assume that the customers who arrive during a glue period of
station $i$ are served first and then customers who retry are
served.

Let $\bar{W}_i$ and $\tilde{W}_i$ denote the waiting time of
type-$i$ customers who arrive during a glue period of station $i$
and any other period, respectively. Further, $G_{i_{res}}$ denotes
the residual time of a glue period of station $i$. Finally
$C_{i_{res}}$ denotes the residual time of a non-glue period of
station $i$. A type-$i$ customer arriving during a glue period of
station $i$ has to wait for the residual glue period. Further, it
has to wait for all the customers who arrived before it during the
glue period. Therefore
\begin{eqnarray*}
    \E [\bar{W_i}]= \E [G_{i_{res}}] + \rho_i \E [G_{i_{res}}] = (1+\rho_i) \E [G_{i_{res}}].
\end{eqnarray*}
A type-$i$ customer arriving during a non-glue period of station
$i$ has to wait for the residual non-glue period, and the glue
period. Then it either gets in the queue for service or it remains
in the orbit. With probability $\tilde G_i (\nu_i)$ it remains in
the orbit and has to wait until the next visit to get served, and
this repeats. Hence, on average, it has to wait for $\tilde G_i
(\nu_i)/(1-\tilde G_i (\nu_i))$ cycles before it gets into the
queue for service. When it gets in the queue it has to wait for
all the type-$i$ customers who have arrived during the glue period
to get served, and then the customers who arrived before it and
who will be served in the current visit period (on average this
number is approximately equal to the number of customers who
arrived during the residual non-glue period before the arrival of
the tagged customer). Therefore
\begin{align*}
    \E [\tilde{W_i}] &\approx \E [C_{i_{res}}] + \E [G_i] + \frac{\tilde G_i (\nu_i)}{1-\tilde G_i (\nu_i)} \E [C]
    + \rho_i \E [G_i] + \rho_i \E[C_{i_{res}}] \\
    &= (1 + \rho_i) \left(\E [C_{i_{res}}] + \E [G_i]\right) + \frac{\tilde G_i (\nu_i)}{1-\tilde G_i (\nu_i)} \E [C].
\end{align*}
The probability that a type-$i$ customer arrives during a glue
period of station $i$ is $\E[G_i]/\mathbb{E}[C]$, and the
probability that it arrives during a non-glue period equals
$1-(\E[G_i]/ \E[C])$. Therefore
\begin{align*}
    \E [W_i] &= \frac{\E[G_i]}{\mathbb{E}[C]} \E [\bar{W_i}] + \frac{\E[C] - \E[G_i]}{\E[C]}\E [\tilde{W_i}] \\
    &\approx (1+\rho_i) \Bigg(\frac{\E[G_i]}{\mathbb{E}[C]} \E [G_{i_{res}}] + \frac{\E[C] - \E[G_i]}{\E[C]}
    \left(\E [C_{i_{res}}] + \E [G_i]\right)\Bigg)+\frac{\tilde G_i (\nu_i)}{1-\tilde G_i (\nu_i)} \left(\E [C] - \E [G_i] \right).
\end{align*}

Let $R_{c_i}$ be the residual cycle time of the system with
respect to station $i$. Then
$$\E[R_{c_i}] = \frac{\E[G_i]}{\mathbb{E}[C]} \E [G_{i_{res}}] + \frac{\E[C] - \E[G_i]}{\E[C]} \left(\E [C_{i_{res}}] + \E [G_i]\right),
\quad i=1,\cdots,N.$$
 We assume that $\E[R_{c_i}] =\E[R_{c}]$ for all $i=1,\ldots,N$. We thus obtain
 (\ref{approx_time_}):
\begin{equation*}
\E [W_i] \approx (1+\rho_i) \E[R_c] + \frac{\tilde G_i
(\nu_i)}{1-\tilde G_i (\nu_i)} \left(\E [C] - \E [G_i]\right).
\end{equation*}

\section*{Acknowledgment}
The research is supported by the IAP program BESTCOM, funded by
the Belgian government, and by the Gravity program NETWORKS,
funded by the Dutch government. The authors gratefully acknowledge
several discussions with Professor Ton Koonen (TU Eindhoven) about
optical communications. B. Kim's research was supported by the
National Research Foundation of Korea (NRF) grant funded by the
Korea government (MSIP) (No. 2014R1A2A2A01005831). J. Kim's
research was supported by Basic Science Research Program through
the National Research Foundation of Korea (NRF) funded by the
Ministry of Education (2014R1A1A4A01003813).

\end{document}